\documentclass[reqno,a4paper,oneside]{amsart}

\usepackage{CJKutf8}
\usepackage{xcolor}
\usepackage{graphicx}
\usepackage[textwidth = 430 pt, textheight = 630 pt]{geometry}

\definecolor{MyDarkBlue}{rgb}{0.15,0.25,0.45}
\definecolor{Cred}{rgb}{0,0,0}
\usepackage{epsfig,rotating}
\usepackage{amsmath,amssymb}
\usepackage{amsfonts}
\usepackage{mathrsfs}
\usepackage{bbm}
\usepackage[normalem]{ulem}
\usepackage{mleftright} 
\usepackage[low-sup]{subdepth} 
\usepackage{latexsym}
\usepackage{amsthm}
\usepackage{tikz}
\usetikzlibrary{matrix,cd,arrows}
\usepackage{lmodern}
\usepackage{microtype}
\usepackage{enumitem} 
\usepackage{setspace} 

\usepackage[all,knot]{xy}
\xyoption{arc}

\usepackage{tikz}
\usetikzlibrary{graphs,decorations.pathmorphing,decorations.markings}
\usepackage{mathtools}
\usepackage[all,knot]{xy}
\xyoption{arc}

\usepackage[toc,page]{appendix}

\renewcommand{\appendices}{
    \section*{Appendix}\label{appendices}\setcounter{subsection}{0}
    \addcontentsline{toc}{section}{Appendix}
    \setcounter{equation}{0}
    \crefalias{subsection}{appendix}
    \makeatletter
    \renewcommand{\theequation}{\Alph{subsection}.\arabic{equation}}
    \renewcommand{\thesubsection}{\Alph{subsection}}
    \renewcommand{\thethm}{\Alph{subsection}.\arabic{thm}}
    \@addtoreset{equation}{subsection}
    \@addtoreset{thm}{subsection}
    \makeatother
}
\let\oldref\ref
\AtBeginDocument{\renewcommand{\ref}[1]{\cref{#1}}}
\AtBeginDocument{\renewcommand{\eqref}[1]{{\rm (\oldref{#1})}}}

\newcommand{\makecommand}[3]{%
    \foreach \i in #3 {%
        \expandafter\xdef\csname #1\i\endcsname{\noexpand#2{\unexpanded\expandafter{\i}}}%
    }%
}

\newcommand{\latinalphabet}{A,a,B,b,C,c,d,D,E,e,F,f,G,g,H,h,I,i,J,j,K,k,L,l,M,m,N,n,O,o,P,p,Q,q,R,r,S,s,T,t,U,u,V,v,W,w,X,x,Y,y,Z,z}

\makecommand{fr}{\mathfrak}{\latinalphabet}

\makecommand{fr}{\mathfrak}{{der,gl,sl,so,osp,brst,su,sp,spin,string,Poisson,inn,at}}

\makecommand{sf}{\mathsf}{\latinalphabet}

\makecommand{sf}{\mathsf}{{id,String,Lie,Hom,SU,inv,SO,Map,Diff,HSymp,id,inn,Inn,CE,hom,Gpd,Bibun,pr}}

\makecommand{rm}{\mathrm}{\latinalphabet}

\makecommand{rm}{\mathrm}{{Ham,ad,id,dim,Ad,im,deg}}

\makecommand{ca}{\mathcal}{\latinalphabet}

\makecommand{I}{\mathbbm}{{N,Z,R,C}}

\makecommand{sc}{\mathscr}{\latinalphabet}

\makecommand{tt}{\mathtt}{\latinalphabet}

\newcommand{\basicc}{\nabla^{\rm{bas}}}

\newcommand{\basiccu}{R_{\nabla}^{\rm{bas}}}

\usepackage{hyperref}
\hypersetup{
    hypertexnames=false,
    colorlinks=true,
    citecolor=MyDarkBlue,
    linkcolor=MyDarkBlue,
    urlcolor=MyDarkBlue,
    pdfauthor={Simon-Raphael Fischer},
    pdftitle={Extension of algebroids I},
    pdfsubject={math.DG},
    breaklinks=true
}

\allowdisplaybreaks

\usepackage[many]{tcolorbox} 
\tcbuselibrary{theorems} 
\usepackage[capitalise,nameinlink,noabbrev]{cleveref} 

\theoremstyle{plain}
\newtheorem{theorem}{Theorem}[section]

\newtcbtheorem
  [use counter*=theorem,number within=section,crefname={lemma}{lemmata},crefname={Lemma}{lemmata}]
  {lemmata}
  {Lemma}
  {%
		fontupper=\itshape,
		breakable,
		enhanced,
		sidebyside adapt=both,
		attach boxed title to top center={yshift=-2mm},
    colback=gray!25,
    colframe=gray!70!black,
    fonttitle=\bfseries,
		subtitle style={boxrule=0.3pt,
		colback=gray!25,
		colupper=black,
		fontupper=\normalfont\footnotesize
		}
  }
  {lem}

\newtcbtheorem
  [use counter*=theorem,number within=section,crefname={proposition}{propositions},crefname={Proposition}{propositions}]
  {propositions}
  {Proposition}
  {%
		fontupper=\itshape,
		breakable,
		enhanced,
		sidebyside adapt=both,
		attach boxed title to top center={yshift=-2mm},
    colback=gray!25,
    colframe=gray!70!black,
    fonttitle=\bfseries,
		subtitle style={boxrule=0.3pt,
		colback=gray!25,
		colupper=black,
		fontupper=\normalfont\footnotesize
		}
  }
  {prop}
\newtcbtheorem
  [use counter*=theorem,number within=section,crefname={theorem}{theorems},crefname={Theorem}{theorems}]
  {theorems}
  {Theorem}
  {%
		fontupper=\itshape,
		breakable,
		enhanced,
		sidebyside adapt=both,
		attach boxed title to top center={yshift=-2mm},
    colback=gray!25,
    colframe=gray!70!black,
    fonttitle=\bfseries,
		subtitle style={boxrule=0.3pt,
		colback=gray!25,
		colupper=black,
		fontupper=\normalfont\footnotesize,
		}
  }
  {thm}
\newtcbtheorem
  [use counter*=theorem,number within=section,crefname={corollary}{corollaries},crefname={Corollary}{corollaries}]
  {corollaries}
  {Corollary}
  {%
		fontupper=\itshape,
		breakable,
		enhanced,
		sidebyside adapt=both,
		attach boxed title to top center={yshift=-2mm},
    colback=gray!25,
    colframe=gray!70!black,
    fonttitle=\bfseries,
		subtitle style={boxrule=0.3pt,
		colback=gray!25,
		colupper=black,
		fontupper=\normalfont\footnotesize
		}
  }
  {cor}

\newtcbtheorem
  [use counter*=theorem,number within=section,crefname={conjecture}{conjectures},crefname={Conjecture}{conjectures}]
  {conjectures}
  {Conjecture}
  {%
		fontupper=\itshape,
		breakable,
		enhanced,
		sidebyside adapt=both,
		attach boxed title to top center={yshift=-2mm},
    colback=gray!25,
    colframe=gray!70!black,
    fonttitle=\bfseries,
		subtitle style={boxrule=0.3pt,
		colback=gray!25,
		colupper=black,
		fontupper=\normalfont\footnotesize
		}
  }
  {con}

\theoremstyle{remark}
\newtcbtheorem
  [use counter*=theorem,number within=section,crefname={remark}{remarks},crefname={Remark}{remarks}]
  {remarks}
  {Remark}
  {%
		breakable,
		enhanced,
		sidebyside adapt=both,
		attach boxed title to top center={yshift=-2mm},
    colback=gray!25,
    colframe=gray!70!black,
    fonttitle=\bfseries,
		subtitle style={boxrule=0.3pt,
		colback=gray!25,
		colupper=black,
		fontupper=\normalfont\footnotesize
		}
  }
  {rem}

\newtheorem{remark}[theorem]{Remarks}

\theoremstyle{definition}

\newtcbtheorem
  [use counter*=theorem,number within=section,crefname={definition}{definitions},crefname={Definition}{definitions}]
  {definitions}
  {Definition}
  {%
		breakable,
		enhanced,
		sidebyside adapt=both,
		attach boxed title to top center={yshift=-2mm},
    colback=gray!25,
    colframe=gray!70!black,
    fonttitle=\bfseries,
		subtitle style={boxrule=0.3pt,
		colback=gray!25,
		colupper=black,
		fontupper=\normalfont\footnotesize
		}
  }
  {def}
\newtcbtheorem
  [use counter*=theorem,number within=section,crefname={example}{examples},crefname={Example}{examples}]
  {examples}
  {Example}
  {
		breakable,
		enhanced,
		sidebyside adapt=both,
		attach boxed title to top center={yshift=-2mm},
    colback=gray!25,
    colframe=gray!70!black,
    fonttitle=\bfseries,
		subtitle style={boxrule=0.3pt,
		colback=gray!25,
		colupper=black,
		fontupper=\normalfont\footnotesize
		}
  }
  {ex}

\makeatletter
\renewcommand{\tocsection}[3]{%
  \indentlabel{\@ifnotempty{#2}{\bfseries\ignorespaces#1 #2\quad}}\bfseries#3}
\renewcommand{\tocsubsection}[3]{%
  \indentlabel{\@ifnotempty{#2}{\ignorespaces#1 #2\quad}}#3}

  \def\l@subsection{\@tocline{2}{0pt}{2.5pc}{5pc}{}}
	
	\makeatother

\begin{document}

		\pagenumbering{Roman}
    \begin{titlepage}
        \vspace*{2cm}
        \begin{center}
            {{\LARGE \bf
                Extension of algebroids Part I: The Construction}
            }
            \vskip1cm
            {
						{\Large {Simon-Raphael \textsc{Fischer}\textsuperscript*}}
						\vskip0.5cm
						1 June 2025: first upload\\
						{\color{Cred}24 March 2026: adding statements regarding crossed modules}
						}
            \vskip0.5cm
            {
                \textsuperscript* Mathematical Institute, Georg-August University Göttingen\\
                Bunsenstr. 3-5, 37073 Göttingen, Germany\\[0.3cm]
            {\href{mailto:simon.fischer@mathematik.uni-goettingen.de}{\texttt{simon.fischer@mathematik.uni-goettingen.de}}}
                }
        \end{center}
        \vskip0.8cm
        \begin{center}
            \textbf{Abstract}\footnote{Abbreviations used in this paper: \textbf{LAB(s)} for Lie algebra bundle(s), \textbf{BLA(s)} for bundle(s) of Lie algebras whose Lie algebras in the fibres may not be isomorphic, \textbf{YM} for Yang-Mills.}
        \end{center}
        \begin{quote}
				In this series of two papers we will generalise the concept of extending a Lie algebroid by a Lie algebra bundle, leading to a notion of extending a Lie algebroid by another Lie algebroid whose orbits lie in the orbits of the former algebroid. The resulting Lie algebroid's anchor will be the sum of the two initial anchors such that the constructions will be similar to matched pairs of Lie algebroids, but with the {\color{Cred} key} difference that we will allow {\color{Cred}a} curvature. In this part of this series we will focus on the canonical construction making use of strict covariant adjustments, a generalisation of Maurer-Cartan forms in the context of gauge theories equipped with a Lie groupoid action instead of a Lie group action. That is, a Cartan connection with certain conditions on the curvature. The second paper will introduce and explain the obstruction of the extension provided here. Examples will include locally split structures as in Poisson geometry {\color{Cred} and crossed modules of Lie groupoids including an obstruction for certain crossed modules on a manifold $M$, generalising the obstruction of a Lie group structure on a manifold}.
				
				As a side result we {\color{Cred} provide an} obstruction theory for certain Cartan connections on Lie algebroids, which will be related to the obstruction of {\color{Cred}what we will call twisted} action algebroids; generalising the statement of the action algebroid structure induced by flat Cartan connections.
        \end{quote}
    \end{titlepage}
		
		 \thispagestyle{empty}
		\hspace{0pt}
		\vfill
		\begin{center}
		\textit{\textbf{This work is dedicated to Kirill C.H.\ Mackenzie. We never met, but without your work I would not be where I am right now, I might not even have been able to finish my Ph.D. Your studies on extending Lie algebroids by Lie algebra bundles helped me understand curved Yang-Mills gauge theories and publishing my first results. It is now my turn to return the favour.
		\\
		\phantom{text}
		\\
		R.I.P.
		}}
		\end{center}
		\vfill
		\newpage
    
		\pagenumbering{arabic}

    \tableofcontents

		\section{Introduction and results}
		
    
		This series of two papers focuses on the extension of a Lie algebroid $F$ by another Lie algebroid $E$ over the same smooth manifold $M$, where one assumes the existence of a Lie algebroid morphism $K \colon E \to F$, in particular, the orbits of the anchor of $E$ are inside the orbits of $F$. The resulting algebroid structure on the {\color{Cred}direct} sum $F \oplus E$ features the sum of anchors on $F$ and $E$ as anchor. This paper focuses on highlighting the construction of the algebroid structure on $F \oplus E$, while the second one will discuss the obstruction and topological invariant behind all of this.
		
		{\color{Cred}Cartan connections (\cite{Abad:0901.0319, Crainic:1210.2277, Blaom:0404313, Abad:0911.2859, Behrend:0410255, Blaom:0509071, Blaom:1605.04365, Crainic:1307.7979, Tang:0405378}) will play a major role in this series of papers, and we will aim to generalise notions related to this type of connection. For this purpose let us quickly summarise some major statements known regarding Cartan connections which we intend to generalise. In this work the abbreviation BLA means ``bundle of Lie algebras,'' that is, a vector bundle $E \to M$ whose fibres are equipped with a Lie bracket giving rise to a tensor field of Lie brackets which we denote by $[\cdot, \cdot]_E$. If the Lie algebras in the fibres are all isomorphic to each other, then we say that $E$ is an LAB, a ``Lie algebra bundle.'' 
		
		Now, if $E$ is a bundle of Lie algebras (BLA), then a Cartan connection on $E$ is a vector bundle connection $\nabla$ such that it is {\color{Cred}a} bracket derivation of $[\cdot, \cdot]_E$; that is, $\nabla[\cdot, \cdot]_E = 0$. This can be canonically generalised to $F$-connections, where $F \to M$ is a Lie algebroid: A Cartan $F$-connection is an $F$-connection $\nabla$ with $\nabla[\cdot, \cdot]_E = 0$. We will keep speaking of a (vector bundle) connection instead of a $\rmT M$-connection, if $F = \rmT M$ is the tangent algebroid of the base manifold.
		
		In this case the following two major milestones are known: 
		\begin{enumerate}[align=left, leftmargin=*]
			\item[BLA i)] A BLA $E$ over a connected manifold $M$ is an LAB if and only if it admits a Cartan connection. (\cite{0521499283, Abad:0901.0319})
			\item[BLA ii)] There is a short exact sequence of Lie algebroids
				\begin{equation*}
					\begin{tikzcd}
						E \arrow[r, hook]& A \arrow[r, two heads] & F ~ ,
					\end{tikzcd}
				\end{equation*}
			if and only if there is a Cartan $F$-connection on $E$ and a $\zeta \in \Gamma\mleft( \wedge^2 F^* \otimes E\mright)$ satisfying 
				\begin{align*}
					R_\nabla &= \mathrm{ad}_E \circ \zeta~,\\
					\mathrm{d}^\nabla \zeta &= 0~,
				\end{align*}
			where $A \to M$ is another Lie algebroid, $R_\nabla$ is the curvature of $\nabla$, ${\rm ad}_E$ is the adjoint representation of $E$, and ${\rmd}^\nabla$ is the exterior covariant derivative of $\nabla$. If $E$ is an LAB, then this is obstructed by what one calls ``Mackenzie's obstruction class.'' (\cite{0521499283}; see also the recent preprint \cite{grad2025covariant} where similar constructions appeared)
		\end{enumerate}
		
		Of course the first statement implies that a BLA $E$ being the kernel of a short exact sequence of Lie algebroids is an LAB along the orbits of the anchor of $F$. 
		
		This series of papers aims to generalise those statements to $E$ being a Lie algebroid. However, one generalisation in that regard is well-known: A Cartan connection on a Lie algebroid $E$ is a vector bundle connection $\nabla$ satisfying
		\begin{equation*}
		0 
		=
		\nabla_X\bigl(\mleft[\mu,\nu\mright]_E\bigr)
			-\mleft[\nabla_X\mu,\nu\mright]_E
			-\mleft[\mu,\nabla_X\nu\mright]_E
			-\nabla_{\basicc_\nu X}\mu
			+\nabla_{\basicc_\mu X}\nu
		\end{equation*}
		for all $\mu ,\nu \in \Gamma(E)$ and $X \in \frX(M)$ (vector field on $M$),
		where $\basicc$ is an $E$-connection on ${\rmT}M$ defined by
		\begin{equation*}
		\basicc_\mu X
		\coloneqq
		\mleft[ \rho_E(\mu), X \mright]_F
			+ \rho_E\mleft( \nabla_X \mu \mright)
		\end{equation*}
		with $\rho_E$ being the anchor of $E$. For completeness, since we will see that the following will play a major role: There is an $E$-connection on $E$ itself also denoted by $\basicc$, 
		\begin{equation*}
		\basicc_\mu \nu
		\coloneqq
		\mleft[ \mu, \nu \mright]_E
			+ \nabla_{\rho_E(\nu)} \mu~,
		\end{equation*}
		such that $\rho_E \circ \basicc = \basicc \circ \rho_E$.
		Using this notion, the following can be derived, a flat algebroid analogue of \textbf{BLA i)} above:
		\begin{enumerate}[align=left, leftmargin=*]
			\item[Flat algebroid i)] A Lie algebroid $E$ over a connected and simply connected manifold $M$ is an action algebroid $M \times \frg$ for some Lie algebra $\frg$ acting on $M$ if and only if $E$ admits a \textbf{flat} Cartan connection (\cite[Thm.~A]{Blaom:0404313}); see also \cite[Prop.\ 2.12]{Abad:0901.0319}, and \cite[Cor.~3.12]{Crainic:1210.2277} for an integrated version.
		\end{enumerate}
		This is a generalisation of a well-known classical statement: A compact, connected and simply connected manifold $M$ admits a Lie group structure if and only if $\rmT M$ admits a flat Cartan connection; see for example \cite[\S 3, Theorem 8.7]{SharpeLieGroupStructure}. The idea is essentially that the corresponding isomorphism $\rmT M \cong M \times \frg$ followed by the canonical projection to $\frg$ is a Maurer-Cartan form, that is, a $\frg$-valued 1-form satisfying the Maurer-Cartan equation which is inherited from the flatness of the Cartan connection.
		
		This will now be the starting point of this series of papers: Observe that the previous statement requires flatness, while the statements regarding BLAs partially include curvature conditions, but not necessarily flatness. This begs the question: 
		
		\begin{center}
		\textit{If flatness is already not required for the BLA statements, what is a proper generalisation of these in the case of Lie algebroids?}
		\end{center}
		
		The tools readily come from what one calls curved Yang-Mills-Higgs theories: Those are gauge theories where a Lie groupoid instead of a Lie group acts on the principal bundle, see \cite{Fischer:2024vak, Fischer:2022sus}; for infinitesimal versions see \cite{Gruetzmann:2014ica, Kotov:2015nuz, Fischer:2020lri, Fischer:2021glc}, and see also \cite{Strobl:2004im, Mayer:2009wf}. In such gauge theories the Maurer-Cartan form and its Maurer-Cartan equation had to be generalised to provide a suitable notion of gauge-invariant theory; that is, there is a theory coming with a generalisation of the Maurer-Cartan form whose curvature condition is not requiring flatness. In addition to this, as pointed out in \cite{Fischer:2020lri, Fischer:2021glc}, whether such gauge theories are ``new'' gauge theories once the Lie algebroid is a BLA is strongly linked to what we described in \textbf{BLA ii)}. Thus there is hope to generalise \textbf{BLA i)} and \textbf{BLA ii)} to Lie algebroids without requiring flatness as Maurer-Cartan equation as in \textbf{Flat algebroid i)}. The generalisation of Maurer-Cartan forms where introduced as \textit{(strict) adjustments} in \cite{Fischer:2024vak} (for the general case) and as \textit{(strict) multiplicative Yang-Mills connections} in \cite{Fischer:2022sus} (for the special case when the structural Lie groupoid is a Lie group bundle).
		
		Indeed, in this series of papers we will investigate strict $K$-adjustment. We will generalise \textbf{BLA i)} and \textbf{Flat algebroid i)}, including a direct generalisation of the classical statement regarding a Lie group structure on a manifold; while we will also already provide the normal form of what will generalise \textbf{BLA ii)}: Given a morphism of Lie algebroids $K\colon E \to F$, for $E, F$ two Lie algebroids over $M$, we will generalise the notion of connection appearing in \textbf{BLA ii)}, by making use of the known definition of Cartan connections on Lie algebroids, as described earlier, while we replace $\rho_E$ with $K$. We will call those connections \textit{Cartan $K$-connections}. In the case of $F = \rmT M$ as the tangent algebroid and $K = \rho_E$ this means that one has an (infinitesimal) multiplicative connection; the curvature $R_\nabla$ is then also multiplicative, but we will require that $R_\nabla$ is multiplicatively exact, giving rise to a $\zeta \in \Gamma\mleft( \wedge^2 F^* \otimes E\mright)$, generalising the Maurer-Cartan equation and the corresponding curvature equation in \textbf{BLA ii)}; this curvature condition can be extended to general $K$ in a straightforward manner. The ``Bianchi identity'' (or strictness) of \textbf{BLA ii)}, $\rmd^\nabla \zeta = 0$, requires a bit more care, but also here, curved Yang-Mills-Higgs theories provide help: By studying the corresponding Atiyah sequence and the obstruction of its algebroid structure. For this we will define another $F$-connection $\nabla^\zeta$ on $E$ given by
		\begin{equation*}
		\nabla^\zeta_X \mu \coloneqq \nabla_X \mu - \zeta \bigl( X, K(\mu) \bigr)
		\end{equation*}
		for all $X \in \Gamma(F)$ and $\mu \in \Gamma(E)$. We will then speak of an \textit{(infinitesimal) strict $K$-adjustment $(\nabla, \zeta)$} if we additionally have
		\begin{equation*}
		\rmd^{\nabla^\zeta} \zeta = 0~.
		\end{equation*}
		If $K = 0$, we will speak of an \textit{(infinitesimal) strict multiplicative Yang-Mills (YM) $F$-connection $(\nabla, \zeta)$} which satisfies the properties listed in \textbf{BLA ii)}; observe that $K = 0$ implies that $E$ is a BLA.
		
		\begin{remarks}{Other references}{OtherLiterature}
		\color{Cred} BLAs and LABs are a special case of Lie algebroids, and, as the previous different labelling of the connections suggests, there is also other literature covering those: On one hand the notion of infinitesimal multiplicative YM connections were of course also introduced in \cite{0521499283}; without the strictness condition they were called ``Lie derivation laws covering a coupling.''
		
		On the other hand there is literature not mentioned above: In \cite{Fischer:2401.05966} strict multiplicative YM connections appeared as an essential tool to classify singular foliations; its appendix classified such connections in the centerless case.
		
		Similar constructions on the group structure also appear in the studies of higher gauge theories, see \cite{Saemann:2019dsl, Kim:2019owc, Rist:2022hci, perez2025higher}.
		
		Multiplicative connections on Lie group or algebra bundles without the curvature conditions also appeared in \cite{Laurent-Gengoux:2005wxa, Castrillon:2201.07088, blazquez2022group, Fernandes:2204.08507}.
		
		Cartan connections on groupoids also appeared in \cite{Chatterjee:2502.02284}, but in a context without curvature conditions; however, the treatment was ``up to homotopy,'' such that \cite{Chatterjee:2502.02284} may actually allow for further generalisations of the results presented here. This may be worked out in future works.
		\end{remarks}
		
		Some major statements are an extension and coordinate-free version of \cite[last statement of Proposition 4.11]{Fischer:2024vak}; the first theorem of this paper is a generalisation of \textbf{BLA i)} and \textbf{Flat algebroid i)}:
		
		\begin{theorems*}{$K$-twisted action algebroid}
		\color{Cred} Let $E$ be a Lie algebroid. First assume that $E$ admits a strict $K$-adjustment $(\nabla, \zeta)$. Then the Lie bracket of $E$ can be written as 
		\begin{equation*}
		\mleft[ \mu, \nu \mright]_E
		=
		H(\mu, \nu)
			+ \nabla^\zeta_{K(\mu)} \nu
			- \nabla^\zeta_{K(\nu)} \mu
			+ \zeta \mleft( K(\mu), K(\nu) \mright) ~ ,
		\end{equation*}
		where $H$ is a field of Lie brackets on $E$ giving rise to a BLA structure given by 
		\begin{equation*}
		H(\mu, \nu)
		=
		t_{\basicc}(\mu, \nu)
			+ \zeta \mleft( K(\mu), K(\nu) \mright)
		\end{equation*}
		for all $\mu, \nu \in \Gamma(E)$, where $t_{\basicc}$ is the torsion of the basic connection (on $E$). Furthermore, $(\nabla^\zeta, \zeta)$ is a strict multiplicative Yang-Mills $F$-connection w.r.t.\ this BLA structure on $E$.
		\newline
		
		Vice versa, given a strict multiplicative Yang-Mills $F$-connection $(\nabla^\zeta, \zeta)$ w.r.t.\ a BLA structure $H$ on $E$ such that the Lie algebroid bracket on $E$ can be written as above, then $(\nabla, \zeta)$ is a strict $K$-adjustment on $E$, where $\nabla$ is given by $\nabla_X \mu \coloneqq \nabla^\zeta_X \mu + \zeta\mleft( X, K(\mu) \mright)$ for all $\mu \in \Gamma(E)$ and $X \in \Gamma(F)$. Furthermore, $H$ can also be written as above.
		\end{theorems*}
		
		Since this will generalise \textbf{Flat algebroid i)} we will call such a structure \textbf{$K$-twisted action algebroid (induced by a strict multiplicative YM connection)} so that this statement can be roughly written as ``There is the structure of a $K$-twisted action algebroid if and only if there is a strict $K$-adjustment'' but one may also want to call it \textbf{twisted crossed module} because we will see in the examples that crossed modules are naturally of this form, too. The BLA structure on $E$ via $H$ will be the \textbf{strict BLA $E_H$ of $E$}, and the duality between $E$ and $E_H$ and between strict adjustments and multiplicative YM connections will simplify finding Cartan connections on algebroids; this may also give hints about how an algebroid $E$ looks like once it comes with a Cartan connection.
		
		Regarding generalising \textbf{BLA ii)} we will derive:
		
		\begin{theorems*}{Sum of algebroids by adjustment}
		\color{Cred} \textit{Let $E, F$ be two Lie algebroids over the same manifold. A strict $K$-adjustment $(\nabla, \zeta)$ on $E$ defines a Lie algebroid structure denoted by $F \bowtie_{(\nabla, \zeta)} E$ on the direct sum $A \coloneqq F \oplus E$ with anchor $\rho_A \coloneqq \rho_F \oplus \rho_E$ and bracket given by
		\begin{align*}
		\mleft[ (X, \mu), (Y, \nu) \mright]_A
		&\coloneqq
		\Bigl( 
			\mleft[ X + K(\mu), Y + K(\nu) \mright]_F
			- K \bigl( 
				\mleft[ \mu, \nu \mright]_E + \nabla_X \nu - \nabla_Y \mu + \zeta(X, Y)
			\bigr)~,
		\\
		&\hspace{1cm}
			\mleft[ \mu, \nu \mright]_E + \nabla_X \nu - \nabla_Y \mu + \zeta(X, Y)
		\Bigr)
		\\
		&=
		\Bigl( 
			\mleft[ X, Y  \mright]_F
			+ \nabla^{\rm bas}_\mu Y
			- \nabla^{\rm bas}_\nu X
			- K \mleft( \zeta(X, Y) \mright) ~,
		\\
		&\hspace{1cm}
			\mleft[ \mu, \nu \mright]_E + \nabla_X \nu - \nabla_Y \mu + \zeta(X, Y)
		\Bigr)
		\end{align*}
		for all $(X, \mu), (Y, \nu) \in \Gamma(A)$. In particular
		\begin{align*}
		\bigl[ (-K(\mu), \mu), (-K(\nu), \nu) \bigr]_A
		&=
		\mleft( 
			- K \mleft( H(\mu, \nu) \mright),
			H(\mu, \nu)
		\mright) ~ ,
		\end{align*}
		where 
		\begin{equation*}
		H(\mu, \nu)
		=
		t_{\basicc}(\mu, \nu)
			+ \zeta \mleft( K(\mu), K(\nu) \mright)
		\end{equation*}
		for all $\mu, \nu \in \Gamma(E)$.}
		\end{theorems*}
		
		We will show that this sits in a sequence which we will call a \textbf{sandglass sequence}:
		\begin{equation*}
			\begin{tikzcd}[column sep=small]
			E \arrow[dd, "K", swap] \arrow[rd, hook] &  & E_H \arrow[ld, swap, hook] \arrow[dd, "-K", dotted]
			\\
			& A \arrow[rd, swap, two heads, dotted] \arrow[dl, two heads] & 
			\\
			F & & F
			\end{tikzcd}
		\end{equation*}
		where arrows with a hook denote injective maps, two heads denote a surjective arrow, and all arrows are Lie algebroid morphisms, except for the dotted ones which are in general only vector bundle morphisms. 
		}
		
		We will conclude this paper with several examples, in particular examples of LABs with multiplicative Yang-Mills connections which appeared already in literature but whose importance may have been neglected. Many examples will come from settings which admit a local splitting theorem like Poisson geometry for example; that is we introduce a natural connection lifting vector fields of a leaf to Poisson vector fields whose curvature will be a Hamilton vector field. These examples will naturally induce a sandglass sequence over the normal bundle of an embedded leaf. {\color{Cred} Last, but not least, we will show that LAB examples naturally induce non-trivial examples in the realm of Lie algebroids: With a suitable LAB action, every Atiyah sequence will naturally induce a (in general) non-trivial sandglass sequence, coming with a strict adjustment on a Lie algebroid which will be given by pulling back the LAB and equipping it with a canonical structure as twisted action Lie algebroid. This will be the major example for strict adjustments.
		
		Other than those examples: Besides the canonical example given by action algebroids we will also highlight that crossed modules of Lie algebroids naturally sit in this framework: Crossed modules are precisely those Lie algebroids which come with a flat adjustment with $\zeta \equiv 0$ and $\basicc \equiv 0$ (both); sandglass sequences will be special types of butterflies in this case, an equivalence of crossed modules as introduced in \cite{Noohi:0910.1818}. Observing this, the last major statement will be a generalisation of the obstruction regarding the existence of a Lie group structure on a compact, connected, and simply connected manifold:
		
		\begin{theorems*}{Obstruction of certain crossed modules of Lie groups}
		\color{Cred}Let $M$ be a smooth, compact, connected and simply connected manifold, $F\to M$ a Lie algebroid, and $K\colon \rmT M \to F$ a morphism of Lie algebroids. 
		
		If $\rmT M$ admits a flat Cartan $K$-connection $\nabla$, then $M$ is diffeomorphic to a Lie group $G$ and there is a crossed module of connected and simply connected Lie groups $(G, \Psi, H, \Phi)$ for which $G$ integrates a Lie algebra induced by a parallel frame of $\basicc$ on $\rmT M$, $H$ integrates a Lie algebra induced by a parallel frame of $\basicc$ on $F$, $K$ and $\nabla$ restrict to these Lie algebras, and $\Psi$ and $\Phi$ are integrals of $K$ and $\nabla$, respectively.

		Furthermore, if we additionally have $\nabla_\xi = 0$ for all $\xi\in \Gamma(F)$ with vanishing anchor, then $H$ sits in a short exact sequence of Lie groups
		\begin{equation*}
			\begin{tikzcd}
				I \arrow[r, hook]& H \arrow[r, two heads] & G ~ ,
			\end{tikzcd}
		\end{equation*}
		where $I$ is the connected and simply connected Lie group integrating the structural isotropy Lie algebra\footnote{As defined by the kernel of $\rho_F$.} of $F$; moreover, $\Psi$ splits this sequence inducing an isomorphism of Lie groups $H \cong G \ltimes I$ (semidirect product).
		\end{theorems*}

		In the second part of this series we will study the obstruction behind sandglass sequences, and this will be linked to the previously mentioned strict adjustments. ``Suitable splittings'' will show that $F\bowtie_{(\nabla, \zeta)} E$ is a normal form of such sequences, and \textit{field redefinitions}, an equivalence relation of curved Yang-Mills-Higgs theories as introduced in \cite{Fischer:2021glc}, will provide the proper notion of change of splittings in this context. In other words: Sandglass sequences will be a generalisation of Mackenzie's studies about short exact sequences of Lie algebroids, motivated by physics.}
		
		The constructions provided here will also strongly resemble the construction of matched pairs of Lie algebroids as in \cite{laurent2008holomorphic} where representations of $F$ on $E$ and vice versa are used, while we only use an $F$-connection $\nabla$ on $E$ which may have a curvature.
		
		\subsection{Notation}
		
		For a Lie algebroid $A$ we denote its anchor by $\rho_A$ and their bracket by $\mleft[ \cdot, \cdot \mright]_A$.
		
		\noindent The base manifold of all involved bundles is usually the smooth manifold $M$, if not otherwise mentioned.
		
		\noindent The sheaf of sections of a bundle $E$ is denoted by $\Gamma(E)$, while vector fields on $M$ are denoted by $\frX(M)$. The sheaf of antisymmetric tensors/forms of degree $k$ are denoted by $\Omega^k$, {\color{Cred}in particular, for two vector bundles $E$ and $F$ over the same base we denote for convenience $\Omega^k(E; F) \coloneqq \Gamma\mleft( \bigwedge^k E^* \otimes F \mright)$, which we call $k$-forms on $E$ with values in $F$. Observe, this is different to the ordinary notation: Forms on $M$ are here denoted by $\Omega^k(\rmT M)$, but this possible confusion will not really play a role here.}
		
		\noindent The action of vector fields $X$ on smooth functions $f$ is denoted by $\scL_X (f)$.
		
		\noindent The total derivative/tangent map of a smooth map $f$ is denoted by $\rmD f$.
		
		
		\subsection{Short exact sequences of algebroids}\label{sec:basics}
		
		Kirill Mackenzie studied short exact sequences of Lie algebroids (\cite{0521499283}), that is, sequences of the form
		\begin{equation*}
			\begin{tikzcd}
				E \arrow[r, hook, "\iota"]& A \arrow[r,"\scD", two heads] & F ~ ,
			\end{tikzcd}
		\end{equation*}
		where $A, E, F$ are Lie algebroids over a smooth manifold $M$, and $\iota, \scD$ are Lie algebroid morphisms. Observe that $E$ is a bundle of Lie algebras (BLA): We have $\scD \circ \iota = 0$, $\rho_F \circ \scD = \rho_A$, and $\rho_A \circ \iota = \rho_E$, thus, altogether, $0 = \rho_F \circ \scD \circ \iota = \rho_A \circ \iota = \rho_E$. This only implies that $E$ is a bundle of Lie algebras, so the fibres might not be isomorphic as Lie algebras; however, {\color{Cred} as mentioned in the introduction, the fibres are isomorphic to each other as Lie algebras along the orbits of $\rho_F$}.

		\begin{remarks}{For now: BLAs}{LABorBLA}
		The mentioned reference assumes that $E$ is an LAB which is needed for the obstruction behind the following. However, we will discuss the obstruction not until the second paper in this series, such that we will not assume that $E$ is an LAB. Observe that the following constructions work for BLAs as well such that it is not required to assume an LAB structure.
		\end{remarks}
		
		According to \cite{0521499283}, a splitting of that sequence is a vector bundle morphism $\chi \colon F \to A$ such that $\scD \circ \chi = {\rm id}_F$, the identity of $F$. It is well-known that such splittings exist and that $\chi$ will be also a morphism of anchored vector bundles; this is simply due to $\rho_A \circ \chi = \rho_F \circ \scD \circ \chi = \rho_F$. However, $\chi$ will be in general not a morphism of Lie algebroids, in particular because its curvature
		\begin{equation*}
		{\color{Cred}\iota}\bigl(\zeta(X, Y) \bigr) \coloneqq {\color{Cred}\iota}\bigl( R_\chi(X, Y) \bigr) \coloneqq 
		\mleft[\chi(X), \chi(Y)\mright]_A
			- \chi\mleft(\mleft[ X, Y \mright]_F\mright)
		\end{equation*}
		is non-zero in general, where $X, Y \in \Gamma(F)$. Since $\scD$ is a morphism of Lie algebroids, {\color{Cred} this is well-defined, and, thus,} $\zeta$ is {\color{Cred}indeed} a 2-form on $F$ with values in $E$, {\color{Cred} i.e.\ $\zeta \in \Omega^2(F; E)$}. Furthermore $\chi$ induces an $F$-connection $\nabla$ on $E$, that is, $\nabla$ is a morphism of anchored vector bundles $F \to \scD(E)$, where $\scD(E)$ is the Lie algebroid of derivations on $E$; for the unfamiliar reader: Such connections behave precisely as typical vector bundle connections except that the Leibniz rule is along $\rho_F$,
		\begin{equation*}
		\nabla_{X} (f\mu) = f \nabla_X \mu + \scL_{\rho_F(X)}(f) ~ \mu
		\end{equation*}
		for all $X \in \Gamma(F)$, $\mu \in \Gamma(E)$, and $f \in C^\infty(M)$; {\color{Cred}we speak just of a ``connection'' if we mean ordinary vector bundle $\rmT M$-connections.} $\nabla$ is defined via
		\begin{equation*}
		{\color{Cred}\iota}\mleft(\nabla_X \mu \mright) \coloneqq \mleft[ \chi(X), {\color{Cred}\iota}(\mu) \mright]_A ~ ,
		\end{equation*}
		which is again well-defined by the fact that $\scD$ is a morphism of Lie algebroids, and that $\chi$ is a morphism of anchored vector bundles. 
		
		\begin{remarks}{Notation}{GenNot}
		{\color{Cred}Denoting $\iota$ is often omitted if the context is clear.}
		\end{remarks}
		
		It is a straightforward exercise to show that 
		\begin{align*}
		\nabla \mleft[ \cdot, \cdot \mright]_E &= 0 ~ ,
		\\
		R_\nabla &= {\rm ad}_E \circ \zeta ~ ,
		\end{align*}
		where ${\rm ad}_E$ denotes the ad-representation in $E$. 
		In fact, one can reconstruct the Lie algebroid structure via $\nabla$ w.r.t.\ the splitting $A \cong F \oplus E$ induced by $\chi$, that is,
		\begin{equation}\label{eq:MackAlgebroid}
		\mleft[ (X, \mu), (Y, \nu) \mright]_A
		=
		\mleft(
			\mleft[X, Y\mright]_F ~ , ~
			\mleft[ \mu, \nu \mright]_E
				+ \nabla_X \nu
				- \nabla_Y \mu
				+ \zeta(X, Y)
		\mright)
		\end{equation}
		for all $(X, \mu), (Y, \nu) \in \Gamma(F \oplus E)$. {\color{Cred} We denote this Lie algebroid structure by $F \ltimes_{(\nabla, \zeta)} E$}, and a change of splitting induces an isomorphism of Lie algebroids between the structures induced by different splittings. For the Jacobi identity it is important that we additionally have $\rmd^\nabla \zeta = 0$, and that is naturally the case for $\zeta = R_\chi$.
		
		In total, \cite{0521499283} shows that $E$ and $F$ sit in such a short exact sequence if and only if there is an $F$-connection $\nabla$ on $E$ with some \textit{primitive} $\zeta$ satisfying
		\begin{align*}
		\nabla \mleft[ \cdot, \cdot \mright]_E &= 0 ~ ,
		\\
		R_\nabla &= {\rm ad}_E \circ \zeta ~ ,
		\\
		\rmd^\nabla \zeta &= 0 ~ .
		\end{align*}
		If $E$ is an LAB, then one can express this compactly in a topological invariant, \textit{Mackenzie's obstruction class} whose role is to measure the existence of a primitive satisfying the third equation once the first two are satisfied. It is then possible to refine the statement to say that ``$E$ and $F$ sit in a short exact sequence covering a coupling if an only if the obstruction is trivial.'' 
		Henceforth, \cite{0521499283} calls $\nabla$ satisfying the first two equations \textit{Lie derivation law covering a coupling\footnote{The term ``coupling'' is mathematically clarified in \cite{0521499283} and sort of provides the existence of $\nabla$ satisfying the first two equations; however we will only reiterate this and the associated statements in the second part of this series since we will not need it here.} between $F$ and $E$}; if the third equation is also satisfied, then we add the adjective \textit{strict}. Before we turn to the agenda of this paper, let us take this moment to come back to $E$ being a BLA in contrast to an LAB. The following is well-known:
		
		\begin{lemmata}{BLA $\stackrel{?}{=}$ LAB}{BLALAB}
		\tcbsubtitle[before skip=\baselineskip]{\cite{0521499283}, see also \cite{Abad:0901.0319}}
		\textit{A BLA $E$ over a connected base manifold is an LAB if and only if it admits a vector bundle connection $\nabla$ such that $\nabla \mleft[ \cdot, \cdot \mright]_E = 0$.}
		\end{lemmata}
		
		That is, if $F = \rmT M$, then the existence of $\nabla$ requires $E$ to be an LAB. However, in general $\nabla$ is an $F$-connection such that we know that the Lie algebra structures of $E$ are the same along the orbits of the anchor of $F$, but the same does not necessarily hold in transversal directions.
		
		{\color{Cred} Since we will not introduce the previously-mentioned couplings here, and to avoid confusion with related terms in physics, we will prefer to call such $\nabla$ as they appear in \cite{Fischer:2022sus, Fischer:2401.05966}; choosing the same label as in the integrated setup:}
		
		\begin{definitions}{Strict multiplicative Yang-Mills (YM) $F$-connections}{StrictYMConn}
		{\color{Cred} Let $E$ be a BLA and $F$ be a Lie algebroid over a smooth manifold $M$. We call a pair $(\nabla, \zeta)$ given by an $F$-connection $\nabla$ on $E$ and a $\zeta \in \Omega^2(F; E)$ a \textit{multiplicative YM $F$-connection}, or just \textit{YM $F$-connection}, if the following are satisfied:}
		\begin{align*}
		\nabla \mleft[ \cdot, \cdot \mright]_E &= 0 ~ ,
		\\
		R_\nabla &= {\rm ad}_E \circ \zeta ~.
		\end{align*}
		We speak of a \textit{strict (multiplicative) YM $F$-connection} if additionally
		\begin{equation*}
			\rmd^\nabla \zeta = 0 ~ .
		\end{equation*}
		\end{definitions}

    \subsection{Strict covariant adjustments}
		
		{\color{Cred}Let us now introduce strict covariant $K$-adjustments, as mentioned in the introduction.} Concretely, we need a special kind of Cartan connection on Lie algebroids, and we will follow\footnote{In the following, certain notions were introduced for $F = \rmT M$ in the mentioned references, but the formulas naturally extend to any Lie algebroid $F$ which is why we will not explicitly recalculate tensorial properties etc.\ here.} closely~\cite{Abad:0901.0319, Crainic:1210.2277, Blaom:0404313}. We work with two Lie algebroids $E \to M$ and $F \to M$ together with an $F$-connection $\nabla$ on $E$. Given a morphism of Lie algebroids $K \colon E \to F$, we define:
		
		\begin{definitions}{Basic connections}{BasicKc}
		\color{Cred} Let $E$ and $F$ be two Lie algebroids over a smooth manifold $M$, and $K\colon E \to F$ a morphism of Lie algebroids. Then we define the \textit{basic $K$-connection $\basicc$ of $\nabla$} as a pair of $E$-connections, one on $E$ itself and the other one on $F$, by
		\begin{align*}
		\basicc_\mu \nu 
		&\coloneqq 
		\mleft[ \mu, \nu \mright]_E
			+ \nabla_{K(\nu)} \mu~,
		\\
		\basicc_\mu X
		&\coloneqq
		\mleft[ K(\mu), X \mright]_F
			+ K\mleft( \nabla_X \mu \mright) ~ ,
		\end{align*}
		\color{Cred}respectively, where $\mu, \nu \in \Gamma(E)$ and $X \in \Gamma(F)$. Usually we will just speak of the basic connection because it will be clear by context what $K$ is.
		\end{definitions}
		
		If $F = \rmT M$, then $K = \rho_E$ is a canonical choice and one recovers the common definition of the basic connection as the infinitesimal version of adjoint (pseudo-)representations in the groupoid setting (upon a choice of $\nabla$). Observe that we have
		\begin{equation*}
		K \circ \basicc = \basicc \circ K ~ .
		\end{equation*}
		In the following we keep the same setup. We also define:
		
		\begin{definitions}{Related notions regarding $\basicc$}{BasicCurvEtc}
		{\color{Cred} The \textit{basic curvature $\basiccu$} is an element of $\Omega^2(E; \operatorname{Hom}(F; E))$ defined by }
		\begin{equation*}
        \basiccu(\mu,\nu)(X)
				\coloneqq 
				\nabla_X\bigl(\mleft[\mu,\nu\mright]_E\bigr)
					-\mleft[\nabla_X\mu,\nu\mright]_E
					-\mleft[\mu,\nabla_X\nu\mright]_E
					-\nabla_{\basicc_\nu X}\mu
					+\nabla_{\basicc_\mu X}\nu~.
    \end{equation*} 
		{\color{Cred}The \textit{torsion} of the basic connection (on $E$) is an element of $\Omega^2(E;E)$ defined by}
		\begin{equation*}
		t_{\basicc}(\mu, \nu) 
		\coloneqq 
		\basicc_\mu \nu 
			- \basicc_\nu \mu 
			- \mleft[ \mu, \nu \mright]_E ~ .
		\end{equation*}
		\end{definitions}
		
		The curvatures of the basic connection on $E$ and on $F$ are equivalent to $- \basiccu \circ K$ and $- K \circ \basiccu$, respectively. 
		and one has (\cite{Blaom:0404313})
		\begin{equation}\label{eq:BasicCurvViaTorsion}
		\basiccu(\mu, \nu)X
		=
		\mleft( \nabla_X t_{\basicc} \mright)(\mu,\nu)
			- R_\nabla\mleft( K(\mu), X \mright) \nu
			+ R_\nabla\mleft( K(\nu), X \mright) \mu ~ ,
		\end{equation}
		where $R_\nabla$ (the curvature of $\nabla$) and terms like $\nabla_X t_{\basicc}$ are defined in the same manner as for vector bundle connections. There is also another canonical $E$-connection on $E$ given by $\nabla_K$, $(\mu, \nu) \mapsto \nabla_{K(\mu)}\nu$, for which one can define the torsion similarly, and it is easy to check that one has 
		\begin{equation*}
		t_{\basicc} = - t_{\nabla_K} ~ .
		\end{equation*}
		
		\begin{definitions}{Cartan $K$-connections}{CartanConn}
		{\color{Cred} We say that $\nabla$ is a \textit{Cartan $K$-connection} if $\basiccu = 0$. If $F= \rmT M$ is the tangent algebroid and $K = \rho_E$, then we just speak of a \textit{Cartan connection}.}
		\end{definitions}
		
		Observe that both $\basicc$ are flat if $\nabla$ is a Cartan $K$-connection. {\color{Cred} Moreover, observe that setting $X = K(\eta)$ in \ref{def:BasicCurvEtc} for $\eta \in E$ implies that $\nabla_K$ is a Cartan ${\rm id}_E$-connection because $\basicc \circ K = K \circ \basicc$ and $(\nabla_K)^{\rm bas} = \basicc$ (on $E$), where one constructs the basic ${\rm id}_E$-connection of $\nabla_K$ here w.r.t.\ ${\rm id}_E$ as Lie algebroid morphism.}
		Let us finally start to state the curvature conditions {\color{Cred}as introduced in \cite[stated for $F= \rmT M$ and $K= \rho_E$]{Fischer:2024vak}.} 
		
		\begin{definitions}{(Strict) covariant $K$-adjustments}{CovAdjust}
		{\color{Cred} A pair $(\nabla, \zeta)$ given by a Cartan $K$-connection $\nabla$ and a $\zeta \in \Omega^2(F; E)$ is a \textit{covariant $K$-adjustment}, or just \textit{$K$-adjustment}, if 
		\begin{equation*}
		R_\nabla = - \rmd^{{\basicc}} \zeta ~ .
		\end{equation*}
		We speak of a \textit{strict (covariant) $K$-adjustment} if additionally
		\begin{equation*}
			\rmd^{\nabla^\zeta} \zeta = 0 ~ ,
		\end{equation*}
		where $\nabla^\zeta$ is an $F$-connection on $E$ defined by
		\begin{equation*}
				\nabla^\zeta_X \nu
				\coloneqq
				\nabla_X \nu
				- \zeta \mleft( X, K(\nu) \mright)
		\end{equation*}
		for all $X \in \Gamma(F)$ and $\nu \in \Gamma(E)$. In both cases we say that $\zeta$ is a \textit{primitive} of $\nabla$, and if $F = \rmT M$ and $K = \rho_E$, then we just write \textit{adjustment} instead of $K$-adjustment.}
		\end{definitions}
		
		\begin{remarks}{Spelling things out}{AdjustmentSpelling}
		{\color{Cred} As usual,} with $\Omega^2(F; E)$ we mean the space of tensors $\Gamma(\bigwedge^2 F^* \otimes E)$. In fact, due to the property that $\basicc$ is a pair of connections, this pair describes an exterior covariant derivative on forms with 2 form-degrees; one w.r.t.\ $F^*$ and another one in $E^*$. In a typical fashion one defines such a covariant derivative on such forms, see \cite[\S 3.8]{Fischer:2021glc}. However, we will not need the general definition, let us just spell out the first curvature condition:
		\begin{equation*}
		R_\nabla(X, Y) \nu 
		= 
		- \basicc_\nu\mleft( \zeta(X, Y) \mright)
			+ \zeta \mleft( \basicc_\nu X, Y \mright)
			+ \zeta \mleft( X, \basicc_\nu Y \mright)
		\end{equation*}
		for all $X, Y \in \Gamma(F)$ and $\nu \in \Gamma(E)$. Since $\zeta$ has no degrees in $E^*$, we can obviously also write $R_\nabla = - \basicc \zeta$; however, we decided against this because we will vary the Cartan $K$-connection in the next paper of this series, and then the notation can be very misleading; see also \cite[in particular Remark 4.5.3]{Fischer:2021glc}. 
		
		{\color{Cred}The exterior covariant derivative $\rmd^{\nabla^\zeta}$ is defined in the same manner as for vector bundle connections. Also, following the comment after \ref{def:CartanConn}, it is straight-forward to check that $\mleft(\nabla_K, \zeta \circ (K, K) \mright)$ is a (strict) ${\rm id}_E$-adjustment if $(\nabla, \zeta)$ is a (strict) $K$-adjustment.}
		\end{remarks}
		
		If $K \equiv 0$, then $E$ has to be a BLA due to $0 = \rho_F \circ K = \rho_E$. In that case {\color{Cred} $\nabla = \nabla^\zeta$}, and the curvature equation reduces to
		\begin{equation*}
		R_\nabla = {\rm ad}_E \circ \zeta ~ ,
		\end{equation*}
		where ${\rm ad}_E$ is the fibre-wise adjoint representation of the BLA $E$; that is, we have a strict YM $F$-connection.
		
		\begin{corollaries}{Strict $0$-adjustments are strict YM $F$-connections}{KEqual0IsYM}
		{\color{Cred} \textit{The notion of (strict) $0$-adjustments is equivalent to the notion of (strict) YM $F$-connections.}}
		\end{corollaries}
		
		{\color{Cred} Last, obviously, a flat Cartan $K$-connection $\nabla$ is a strict $K$-adjustment by choosing $\zeta \equiv 0$ as primitive, and we usually choose such a $\zeta$ if not mentioned otherwise.}
		
		\section{Action algebroids via multiplicative YM \texorpdfstring{$F$}{F}-connections}
		
		{\color{Cred} Let us start with the first theorem of this paper which will generalise the bullet point \textbf{Flat algebroid i)} of the introduction. Henceforth, for readability, let us first revisit the proof of \textbf{Flat algebroid i)}, and let us observe that we actually do not need to assume that $M$ is connected and simply connected.}
		
		\begin{propositions}{Flat case revisited I}{RevisitObstrActAlgoid}
		{\color{Cred}\textit{Let $E, F \to M$ be two Lie algebroids, $K\colon E \to F$ a morphism of Lie algebroids, coming with a flat Cartan $K$-connection $\nabla$. Then the torsion $t_{\basicc}$ defines a BLA structure on $E$ with $\nabla t_{\basicc} = 0$, and we can write
		\begin{equation*}
			\mleft[ \mu, \nu \mright]_E
			=
			t_{\basicc}(\mu, \nu)
				+ \nabla_{K(\mu)} \nu
				- \nabla_{K(\nu)} \mu
		\end{equation*}
		for all $\mu, \nu \in \Gamma(E)$.}}
		\end{propositions}
		
		{\color{Cred}
		\begin{proof} 
		Due to the flatness of $\nabla$, Eq.\ \eqref{eq:BasicCurvViaTorsion} reduces to
		\begin{equation*}
		0 = \nabla t_{\basicc}~.
		\end{equation*}
		The classical proof (for $F = \rmT M$ and $K=\rho_E$) would now restrict to a local situation by fixing a parallel frame, that is, one makes use of the flatness of $\nabla$ to conclude the action algebroid structure via the behaviour of the torsion now. However, we do not need to restrict to a local statement: The famous Bianchi identity relates curvatures and torsions for arbitrary setups, so let us rather make use of this identity. In the same fashion as for vector bundle connections, it is straightforward to check that one has the following Bianchi identity for $\nabla_{K}$ as an $F$-connection on $E$, i.e.\
		\begin{equation*}
		0
		=
		\sum_{\text{cyclic in } \mu, \nu, \eta} \mleft(
			t_{\nabla_{K}}\mleft( t_{\nabla_{K}}(\mu, \nu), \eta \mright)
			+ \mleft( \nabla_{{K}(\mu)} t_{\nabla_{K}} \mright)(\nu, \eta)
			- R_{\nabla_{K}}(\mu, \nu) \eta
		\mright)
		\end{equation*}
		for all $\mu, \nu, \eta \in \Gamma(E)$. $\nabla_{K}$ is flat due to the flatness of $\nabla$ and because $K$ is a morphism of Lie algebroids, and, as previously mentioned, $t_{\basicc} = - t_{\nabla_{K}}$, so that we can combine this with the previous equation to conclude that the Jacobiator of $t_{\basicc}$ is zero,
		\begin{equation*}
		0 = \sum_{\text{cyclic in } \mu, \nu, \eta} 
			t_{\basicc}\mleft( t_{\basicc}(\mu, \nu), \eta \mright)~.
		\end{equation*}
		That is, as a Lie bracket, $t_{\basicc}$ gives rise to a BLA structure on $E$. That $[\cdot, \cdot]_E$ can be written as proposed is just a consequence of the definition of the torsion $t_{\basicc} = - t_{K}$, solved for $[\cdot, \cdot]_E$.
		\end{proof}
		
		For $F = \rmT M$ and $K=\rho_E$ we can therefore immediately conclude:
		
		\begin{corollaries}{Flat case revisited II: Obstruction of trivial action algebroids}{RevisitObstrActAlgoid}
		{\color{Cred}\textit{Let $E \to M$ be a Lie algebroid, equipped with a flat Cartan connection $\nabla$. Then the torsion $t_{\basicc}$ defines an LAB structure on $E$ with $\nabla t_{\basicc} = 0$, and we can write
		\begin{equation*}
			\mleft[ \mu, \nu \mright]_E
			=
			t_{\basicc}(\mu, \nu)
				+ \nabla_{\rho_E(\mu)} \nu
				- \nabla_{\rho_E(\nu)} \mu
		\end{equation*}
		for all $\mu, \nu \in \Gamma(E)$.}}
		\end{corollaries}
		
		\begin{proof}
		Just a consequence of \ref{prop:RevisitObstrActAlgoid}, while the LAB structure follows by \ref{lem:BLALAB}.
		\end{proof}
		
		W.r.t.\ parallel frames this reduces to the typical definition action Lie algebroids defined by Lie algebra actions. Thus, this motivates the first theorem: By allowing a curvature we adjust the torsion in such a way that it is a Lie bracket, while the curvature conditions assure the Jacobi identity; that we can express the Lie bracket of the algebroid structure by this structure will be again just a direct consequence of the involved definitions.
		
		In order to do so in a concise way, let us first state the outcoming Lie algebroid structure:

		\begin{propositions}{Formal definition of $K$-twisted action algebroids}{defTwistActionAlg}
		\color{Cred} \textit{Let $E_H \to M$ be a BLA with field of Lie brackets $H$, $F \to M$ a Lie algebroid, and let $(\nabla, \zeta)$ be a strict YM $F$-connection on $E_H$. Furthermore assume there is a vector bundle morphism $K \colon E_H \to F$, satisfying}
		\begin{equation*}
		\mleft[ K(\mu), K(\nu) \mright]_F
		=
		K \mleft(
			H(\mu, \nu)
			+ \nabla_{K(\mu)} \nu
			- \nabla_{K(\nu)} \mu
			+ \zeta \mleft( K(\mu), K(\nu) \mright)
			\mright)
		\end{equation*}
		\textit{for all $\mu, \nu \in \Gamma(E_H)$. Then $E_H$ admits the structure as a $K$-twisted action algebroid\footnote{Thanks to the referee for the suggested label.} (induced by $H$ and $(\nabla, \zeta)$), which we denote as $E$ and for which $K \colon E \to F$ is a morphism of Lie algebroids, that is, the anchor is given by $\rho_F \circ K$, and the Lie bracket is given as }
		\begin{equation*}
		\mleft[ \mu, \nu \mright]_E
		=
		H(\mu, \nu)
			+ \nabla_{K(\mu)} \nu
			- \nabla_{K(\nu)} \mu
			+ \zeta \mleft( K(\mu), K(\nu) \mright) ~ ,
		\end{equation*}
		\textit{for all $\mu, \nu \in \Gamma(E)$. We speak of a $K$-twisted action BLA (LAB) if $\rho_F \circ K \equiv 0$.}
		\end{propositions}
		
		\begin{examples}{Recovering ``trivial'' action Lie algebroids}{RecoverTrivialActionAlgoid}
		\color{Cred} Observe that we recover the definition of ``trivial'' action algebroids as motivated by \ref{cor:RevisitObstrActAlgoid}: Let $\frg$ be a Lie algebra with a linear map $\gamma\colon \frg \to \frX(M)$, then define:
		\begin{itemize}
			\item $F = \rmT M$, the tangent algebroid of $M$,
			\item $E_H \coloneqq M \times \frg$ the trivial LAB with canonical bracket $H \coloneqq \mleft[ \cdot, \cdot \mright]_\frg$,
			\item \color{Cred}$\zeta \equiv 0$,
			\item $\nabla \coloneqq \nabla^0$, the canonical flat connection so that $(\nabla^0, \zeta\equiv 0)$ is the canonical strict YM connection on $E_H$,
			\item $K(p, \mu) \coloneqq \gamma(\mu)_p$ for all $(p, \mu) \in M \times \frg$.
		\end{itemize}
		\color{Cred}The needed condition on $K$ is equivalent to $\gamma$ being a Lie algebra action: The curvature $R_K$ of $K$, $R_K(\mu, \nu) = \mleft[ K(\mu), K(\nu) \mright]_F - K\mleft(\mleft[ \mu, \nu \mright]_E\mright)$ for $\mu, \nu \in \Gamma(E)$, is a tensor because, by definition, $K$ is a morphism of anchored vector bundles $E \to F$ and $\mleft[ \cdot, \cdot \mright]_E$ comes with the canonical Leibniz rule with $\rho_F \circ K$ as anchor, and thus $R_K$ can be tested at constant sections. Doing so, we see that $R_K = 0$ if and only if $\gamma$ is an action of $\frg$.
		\end{examples}
		
		\begin{remarks}{General action Lie algebroids}{generalActionAlgoidsNotNice}
		\color{Cred} While we recover action Lie algebroids for Lie algebra actions, we do not recover the full notion of action Lie algebroids. A Lie algebroid $E$ over a manifold $N$ acting on another manifold $M$ along some map $\Phi\colon M \to N$ may not be of such a form, for example, if $E$ is a BLA, then so is $\Phi^*E$, but, by \ref{lem:BLALAB}, $\Phi^*E$ may not admit a strict YM connection. Take for example the adjoint action of $E$ on itself, while $E$ is a BLA which is not an LAB.
		
		However we will see in the examples that action Lie algebroids have this form if $E$ itself already comes with a strict adjustment.
		\end{remarks}
		
		\begin{proof}[Proof of \ref{prop:defTwistActionAlg}]
		By assumption, $K$ is a homomorphism of antisymmetric bilinear brackets. Now, in order to show the Jacobi identity one makes use of the proposed expression of the Lie algebroid bracket which resembles the second component of Eq.\ \eqref{eq:MackAlgebroid}: It is a standard exercise to show that the Jacobi identity of $[\cdot, \cdot]_E$ {\color{Cred} follows by} the following set of equations: Jacobi identity of $H$, and the equations\footnote{Along the image of $K$ because the involved contractions in these equations are with sections of the form $K(\mu)$, and not with general sections of $F$.} $\nabla H = 0$, $R_{\nabla} = {\rm ad}_H \circ \zeta$, and $\rmd^{\nabla} \zeta = 0$; {\color{Cred}this is in the same manner as in \cite[proof of Theorem 7.3.7]{0521499283} where the Jacobi identity of Eq.\ \eqref{eq:MackAlgebroid} is proven, that is, following the notation of Eq.\ \eqref{eq:MackAlgebroid}, set $X= K(\mu)$ and $Y=K(\nu)$, the proof of \cite{0521499283} then carries over by making use of the fact that $K$ is a homomorphism of brackets.}
		\end{proof}
		
		}
		
		\begin{theorems}{Duality of adjustments and $K$-twisted action algebroids}{ActionAlgebroidStructure}
		\tcbsubtitle[before skip=\baselineskip]{For the first bullet point: \cite[especially the last statement of Proposition 4.11]{Fischer:2024vak}}
		\color{Cred}
		We have the following duality:
		\begin{itemize}
			\item Let $E, F$ be two Lie algebroids over a smooth manifold $M$, and $K \colon E \to F$ a morphism of Lie algebroids. If there is a strict $K$-adjustment $(\nabla, \zeta)$, then $E$ admits a BLA structure $E_H$ with a field of Lie brackets $H$ given by
			\begin{equation*}
			H(\mu, \nu)
			=
			t_{\basicc}(\mu, \nu)
				+ \zeta \bigl( K(\mu), K(\nu) \bigr)
			\end{equation*}
			for all $\mu, \nu \in \Gamma(E)$, and $(\nabla^\zeta, \zeta)$ is a strict YM $F$-connection on $E_H$. Furthermore, $E$ is then the $K$-twisted action algebroid structure induced by $H$ and $(\nabla^\zeta, \zeta)$.
			
			\item Vice versa, given a $K$-twisted action algebroid $E$ as in \ref{prop:defTwistActionAlg} induced by a strict YM $F$-connection $(\nabla^\zeta, \zeta)$ w.r.t.\ a BLA structure $H$ on $E$, then $(\nabla, \zeta)$ is a strict $K$-adjustment on $E$, where $\nabla$ is given by $\nabla_X \mu \coloneqq \nabla^\zeta_X \mu + \zeta\mleft( X, K(\mu) \mright)$ for all $\mu \in \Gamma(E)$ and $X \in \Gamma(F)$. Furthermore, $H$ can also be written as above.
		\end{itemize}
		\end{theorems}
		
		\begin{proof}[Proof of \ref{thm:ActionAlgebroidStructure}] 
		The first part has already been shown in \cite[a consequence of the last statement of Proposition 4.11; see also the discussion in the next section here]{Fischer:2024vak}, thus let us first focus on the second part. However, we will clarify at this proof's end how the first part can also be shown by reverting the following calculation such that the reference is not needed. 
		
		{\color{Cred}By \ref{prop:defTwistActionAlg} }we have
		\begin{align*}
		\nabla_X \mleft( \mleft[ \mu, \nu \mright]_E \mright)
		&=
		\nabla^\zeta_X \mleft( \mleft[ \mu, \nu \mright]_E \mright)
			 + \zeta \mleft( X, K\mleft( \mleft[ \mu, \nu \mright]_E \mright) \mright)
		\\
		&=
		H\mleft(\nabla^\zeta_X\mu, \nu\mright)
		+ H\mleft(\mu, \nabla^\zeta_X \nu\mright)
		\\
		&\hspace{1cm}
			+ \nabla^\zeta_X \nabla^\zeta_{K(\mu)} \nu
			- \nabla^\zeta_X \nabla^\zeta_{K(\nu)} \mu
		\\
		&\hspace{1cm}
			+ \nabla^\zeta_X \mleft(\zeta \mleft( K(\mu), K(\nu) \mright) \mright)
			+ \zeta \mleft( X, K\mleft( \mleft[ \mu, \nu \mright]_E \mright) \mright) 
		\end{align*}
		{\color{Cred}for all $\mu, \nu \in \Gamma(E)$ and $X \in \Gamma(F)$,}
		where we made use of $\nabla^\zeta H = 0$. 
		We also have
		\begin{align*}
		\mleft[ \nabla_X\mu, \nu \mright]_E
		&=
		H\mleft(\nabla_X\mu, \nu\mright)
			+ \nabla^\zeta_{K\mleft(\nabla_X\mu\mright)} \nu
			- \nabla^\zeta_{K(\nu)} \nabla_X\mu
			+ \zeta \mleft( K\mleft(\nabla_X\mu \mright), K(\nu) \mright)
		\\
		&=
		H\mleft(\nabla^\zeta_X\mu, \nu\mright)
			+ \nabla_{K\mleft(\nabla_X\mu\mright)} \nu
			- \nabla^\zeta_{K(\nu)} \nabla^\zeta_X\mu
		\\
		&\hspace{1cm}
			+ H\mleft(\zeta (X, K(\mu)), \nu\mright)
			- \nabla^\zeta_{K(\nu)} \mleft( \zeta \mleft( X, K(\mu) \mright) \mright) ~ ,
		\end{align*}
		and 
		\begin{equation*}
		\nabla_{\basicc_\mu X} \nu
		=
		\nabla_{\mleft[ K(\mu), X \mright]_F} \nu
			+ \nabla_{K(\nabla_X \mu)} \nu
		=
		\nabla^\zeta_{\mleft[ K(\mu), X \mright]_F} \nu
			+ \nabla_{K(\nabla_X \mu)} \nu
			+ \zeta\mleft( \mleft[ K(\mu), X \mright]_F, K(\nu) \mright) ~.
		\end{equation*}
		Thus,
		\begin{align*}
		\basiccu(\mu,\nu)(X)
		&=
		\nabla_X\mleft(\mleft[\mu,\nu\mright]_E\mright)
			-\mleft[\nabla_X\mu,\nu\mright]_E
			-\mleft[\mu,\nabla_X\nu\mright]_E
			-\nabla_{\basicc_\nu X}\mu
			+\nabla_{\basicc_\mu X}\nu
		\\
		&=
		\nabla^\zeta_X \nabla^\zeta_{K(\mu)} \nu
			- \nabla^\zeta_X \nabla^\zeta_{K(\nu)} \mu
			+ \nabla^\zeta_X \mleft(\zeta \mleft( K(\mu), K(\nu) \mright) \mright)
			+ \zeta \mleft( X, \mleft[ K(\mu), K(\nu) \mright]_F \mright)
		\\
		&\hspace{1cm}
			+ \nabla^\zeta_{K(\nu)} \nabla^\zeta_X\mu
			- H\mleft(\zeta (X, K(\mu)), \nu\mright)
			+ \nabla^\zeta_{K(\nu)} \mleft( \zeta \mleft( X, K(\mu) \mright) \mright)
		\\
		&\hspace{1cm}
			- \nabla^\zeta_{K(\mu)} \nabla^\zeta_X\nu
			+ H\mleft(\zeta (X, K(\nu)), \mu\mright)
			- \nabla^\zeta_{K(\mu)} \mleft( \zeta \mleft( X, K(\nu) \mright) \mright)
		\\
		&\hspace{1cm}
			- \nabla^\zeta_{\mleft[ K(\nu), X \mright]_F} \mu
			- \zeta\mleft( \mleft[ K(\nu), X \mright]_F, K(\mu) \mright)
		\\
		&\hspace{1cm}
			+ \nabla^\zeta_{\mleft[ K(\mu), X \mright]_F} \nu
			+ \zeta\mleft( \mleft[ K(\mu), X \mright]_F, K(\nu) \mright)
		\\
		&=
		\underbrace{R_{\nabla^\zeta}\mleft( X, K(\mu) \mright)\nu
			- H\mleft(\zeta (X, K(\mu)), \nu\mright)}_{= 0}
		\\
		&\hspace{1cm}
			\underbrace{- R_{\nabla^\zeta}\mleft( X, K(\nu) \mright)\mu
			+ H\mleft(\zeta (X, K(\nu)), \mu\mright)}_{= 0}
		\\
		&\hspace{1cm}
			+ \underbrace{\mleft( \rmd^{\nabla^\zeta} \zeta \mright)}_{=0} \mleft( X, K(\mu), K(\nu) \mright)
		\\
		&=
		0 ~ .
		\end{align*}
		Thus, $\nabla$ is a Cartan connection; it is only left to show that it is a covariant adjustment, strictness already follows by assumption. Observe that we can also write
		\begin{align*}
		H(\mu, \nu) 
		&=
		[\mu ,\nu]_E
			- \nabla^\zeta_{K(\mu)} \nu
			+ \nabla^\zeta_{K(\nu)} \mu
			- \zeta \mleft( K(\mu), K(\nu) \mright)
		\\
		&=
		- t_{\nabla_K}(\mu, \nu)
			+ \zeta \mleft( K(\mu), K(\nu) \mright)
			- \zeta \mleft( K(\nu), K(\mu) \mright)
			- \zeta \mleft( K(\mu), K(\nu) \mright)
		\\
		&=
		- t_{\nabla_K}(\mu, \nu)
			+ \zeta \mleft( K(\mu), K(\nu) \mright) ~ ,
		\end{align*}
		such that on $E$
		\begin{equation*}
		\basicc_\mu \nu 
		=
		\nabla_{K(\mu)} \nu 
			- t_{\nabla_K}(\mu, \nu)
		= 
		\nabla_{K(\mu)} \nu 
			+ H(\mu, \nu)
			- \zeta\mleft( K(\mu), K(\nu) \mright)
		=
		\nabla^\zeta_{K(\mu)} \nu 
			+ H(\mu, \nu)~,
		\end{equation*}
		but recall that we have on $F$
		\begin{equation*}
		\basicc_\mu X
		=
		\mleft[ K(\mu), X \mright]_F
			+ K \mleft( \nabla_X \mu \mright) ~,
		\end{equation*}
		thus,
		\begin{align*}
		\nabla_X \nabla_Y \mu
		&=
		\nabla^\zeta_X \nabla^\zeta_Y \mu
			+ \nabla^\zeta_X \mleft( \zeta(Y, K(\mu) \mright)
			+ \zeta \mleft( X, K \mleft( \nabla_Y \mu \mright) \mright)
		\\
		&=
		\nabla^\zeta_X \nabla^\zeta_Y \mu
			+ \nabla^\zeta_X \mleft( \zeta(Y, K(\mu) \mright)
			+ \zeta \mleft( X, \basicc_\mu Y \mright)
			- \zeta \mleft( X, \mleft[ K(\mu), Y \mright]_F \mright) ~.
		\end{align*}
		In total
		\begin{align*}
		R_\nabla(X, Y) \mu
		&=
		R_{\nabla^\zeta}(X, Y) \mu 
			- \zeta \mleft( Y, \basicc_\mu X \mright)
			+ \zeta \mleft( X, \basicc_\mu Y \mright)
		\\
		&\hspace{1cm}
			+ \nabla^\zeta_X \mleft( \zeta(Y, K(\mu) \mright)
			- \nabla^\zeta_Y \mleft( \zeta(X, K(\mu) \mright)
		\\
		&\hspace{1cm}
			- \zeta\mleft( \mleft[ X, Y \mright]_F, K(\mu) \mright)
			+ \zeta \mleft( \mleft[ X, K(\mu) \mright]_F, Y \mright)
			- \zeta \mleft( \mleft[ Y, K(\mu) \mright]_F, X \mright)
		\\
		&=
		- H\mleft( \mu, \zeta(X, Y) \mright)
			+ \zeta \mleft( \basicc_\mu X, Y \mright)
			+ \zeta \mleft( X, \basicc_\mu Y \mright)
		\\
		&\hspace{1cm}
			+ \rmd^{\nabla^\zeta}\zeta\mleft( X, Y, K(\mu) \mright)
			- \nabla^\zeta_{K(\mu)} \mleft( \zeta(X, Y \mright)
		\\
		&=
		- \basicc_\mu\mleft( \zeta(X, Y) \mright)
			+ \zeta \mleft( \basicc_\mu X, Y \mright)
			+ \zeta \mleft( X, \basicc_\mu Y \mright)
		\\
		&=
		-\rmd^{\basicc} \zeta (X, Y, \mu) ~.
		\end{align*}
		{\color{Cred} This proves the second bullet point of \ref{thm:ActionAlgebroidStructure}. Regarding the first bullet point:} Following this calculation there is an alternative proof to \cite[especially the last statement of Proposition 4.11]{Fischer:2024vak}. Just revert the previous calculations but starting with $R_\nabla$ for which we actually have shown that 
		\begin{equation*}
		R_\nabla(X, Y) \mu 
		=
		-\rmd^{\basicc} \zeta (X, Y, \mu) 
			+ R_{\nabla^\zeta}(X, Y) \mu 
			- H\mleft( \zeta(X, Y), \mu \mright) ~ ,
		\end{equation*}
		such that the curvature equation for $\nabla^\zeta$ promptly follows. Reverting the calculations above to show $\nabla^\zeta H = 0$ and that $H$ is a field of Lie brackets is then straightforward: One starts with the same calculation as in this proof, which is possible because the proposed expression for the Lie algebroid bracket {\color{Cred}as provided in \ref{prop:defTwistActionAlg}} is just a consequence of the definitions {\color{Cred}and $H = t_{\basicc} + \zeta \circ (K,K)$}; but now carry the $\nabla^\zeta H$ term along the way. Since strictness and the curvature equation for $\nabla^\zeta$ is at this step known, it follows immediately that $\nabla^\zeta H = 0$ by the fact that $R_\nabla^{\rm bas} = 0$. 
		
		$H$ is by construction antisymmetric and a tensor; in order to show the Jacobi identity one makes again use of the proposed expression of the Lie algebroid bracket. {\color{Cred} That is, the Jacobi identity follows as in the proof of \ref{prop:defTwistActionAlg}, but by reverting the argument, i.e.,} one similarly argues that the Jacobi identity of $H$ follows by the Jacobi identity of $[ \cdot, \cdot ]_E$, $\nabla^\zeta H = 0$, $R_{\nabla^\zeta} = {\rm ad}_H \circ \zeta$, $\rmd^{\nabla^\zeta} \zeta = 0$, and the fact that $K$ is a morphism of Lie algebroids, all of which are either already known or shown.
		\end{proof}
		
		{\color{Cred} This allows us to give a direct generalisation of \textbf{Flat algebroid i)} of the introduction, but observe that the following statement still requires curvature conditions while \textbf{BLA i)} of the introduction does not.}
		
		\begin{corollaries}{Obstruction of $K$-twisted action algebroids}{ObsTwistedActionAlgoid}
		\color{Cred} \textit{Let $E, F$ be two Lie algebroids over a smooth manifold $M$, and $K \colon E \to F$ a morphism of Lie algebroids. Then $E$ is a $K$-twisted action algebroid induced by some BLA structure coming with a strict YM $F$-connection if and only if it admits a strict $K$-adjustment.}
		\end{corollaries}

		
		
		\begin{remarks}{Longitudinal LAB structure}{LABalongOrbits}
		In general, the Lie algebra structures in each fibre of $E$ via $H$ are only isomorphic to each other along the orbits of the anchor of $F$, recall the discussion around \ref{lem:BLALAB}. In particular, if $F$ is transitive (for example $F = \rmT M$), then $E$ admits an LAB structure. This will be important to understand once we turn to the obstruction class behind that.
		\end{remarks}
		
		Let us conclude this section by giving the BLA structure of $E$ a name.
		
		\begin{definitions}{Strict LAB structure}{StrictLABs}
		The BLA (LAB) structure on {\color{Cred}a $K$-twisted action algebroid $E$ is} the \textit{strict BLA (LAB) structure of $E$}. {\color{Cred} As before,} we write $E_H$ instead of $E$ if we want to speak of $E$ as a BLA (LAB), where $H$ is its field of Lie brackets.
		\end{definitions}
		
		\begin{propositions}{$H$ is constant w.r.t.\ $\basicc$}{HConstantWRTthebasicConn}
		{\color{Cred}Let $E, F$ be two Lie algebroids over a smooth manifold $M$, $K \colon E \to F$ a morphism of Lie algebroids, and $(\nabla, \zeta)$ a $K$-adjustment. Then} we have
		\begin{equation*}
		\basicc H = 0~,
		\end{equation*}
		{\color{Cred}where $H(\mu, \nu) = t_{\basicc}(\mu, \nu) + \zeta \bigl( K(\mu), K(\nu) \bigr)$ for all $\mu, \nu \in \Gamma(E)$.}
		\end{propositions}
		
		\begin{proof}
		This is a straightforward consequence of \cite[{\color{Cred}Corollary 3.6.6, Lemma 3.8.5}, Theorem 4.8.4]{Fischer:2021glc}, based on \cite{Blaom:0404313}. That is, by these references, on one hand one has
		\begin{equation*}
		R_{\nabla_K} = \basicc t_{\basicc} ~ ,
		\end{equation*}
		on the other hand
		\begin{equation*}
		R_{\nabla_K}
		=
		\mleft( - \rmd^{\basicc} \zeta \mright) \circ (K, K)
		=
		- \basicc\mleft( \zeta \circ (K, K) \mright) ~ ,
		\end{equation*}
		and thus $\basicc H = 0$; see the references for the involved calculations.
		\end{proof}
		
		\begin{remarks}{Another YM connection and a possible primitive along $K$}{PossibleChoiceOfZetaKK}
		{\color{Cred} In particular, since $\basicc$ is flat, $\basicc$ (on $E$) canonically describes a strict YM $E$-connection on $E_H$. If $\zeta \equiv 0$, then $E$ is trivially an $\mathrm{id}_E$-twisted action algebroid due to $[\mu,\nu]_E = - t_{\basicc}(\mu, \nu) + \basicc_\mu \nu - \basicc_\nu \mu$ for all $\mu, \nu \in \Gamma(E)$; that is, we set the BLA bracket now as $-H = - t_{\basicc} = t_{\nabla_K}$. By \ref{thm:ActionAlgebroidStructure} $(\basicc, \zeta \equiv 0)$ then describes a flat strict $\mathrm{id}_E$-adjustment. This is useful to observe once we discuss crossed modules and \ref{thm:ObstrCrossedModules}.}
		
		As one {\color{Cred}also} sees in that proof, a possible choice for $\zeta$ along the orbits of $K$ may be $- t_{\basicc} = t_K$ {\color{Cred} if one can lift $F$ into $E$ following $K$}, inducing an abelian structure on $E_H$. Indeed, this corresponds to a sort of Yang-Mills-Higgs theory induced by an abelian Lie algebra action; see \cite[Corollary 4.4.9 and Corollary 4.8.5]{Fischer:2021glc}.
		\end{remarks}
		
		\section{Sandglass sequences}\label{sec:Construction}
		
		Finally, we can state the leading construction in this paper {\color{Cred}which will serve as a normal form of what we will call sandglass sequence:}
		
		\begin{theorems}{Sum of algebroids by adjustment}{CanonEx}
		\tcbsubtitle[before skip=\baselineskip]{\cite[especially the last statement of Proposition 4.11]{Fischer:2024vak}}
		\textit{{\color{Cred}Let $E, F$ be two Lie algebroids over a smooth manifold $M$, and $K \colon E \to F$ a morphism of Lie algebroids.} A strict covariant $K$-adjustment $(\nabla, \zeta)$ on $E$ defines a Lie algebroid structure {\color{Cred} denoted as $F \bowtie_{(\nabla, \zeta)}E$} on the {\color{Cred}direct} sum $A \coloneqq F \oplus E$ with anchor $\rho_A \coloneqq \rho_F \oplus \rho_E$ and bracket given by}
		\begin{align*}
		\mleft[ (X, \mu), (Y, \nu) \mright]_A
		&\coloneqq
		\Bigl( 
			\mleft[ X + K(\mu), Y + K(\nu) \mright]_F
			- K \bigl( 
				\mleft[ \mu, \nu \mright]_E + \nabla_X \nu - \nabla_Y \mu + \zeta(X, Y)
			\bigr)~,
		\\
		&\hspace{1cm}
			\mleft[ \mu, \nu \mright]_E + \nabla_X \nu - \nabla_Y \mu + \zeta(X, Y)
		\Bigr)
		\\
		&=
		\Bigl( 
			\mleft[ X, Y  \mright]_F
			+ \nabla^{\rm bas}_\mu Y
			- \nabla^{\rm bas}_\nu X
			- K \mleft( \zeta(X, Y) \mright) ~,
		\\
		&\hspace{1cm}
			\mleft[ \mu, \nu \mright]_E + \nabla_X \nu - \nabla_Y \mu + \zeta(X, Y)
		\Bigr)
		\end{align*}
		\textit{for all $(X, \mu), (Y, \nu) \in \Gamma(A)$. In particular}
		\begin{align*}
		\bigl[ (-K(\mu), \mu), (-K(\nu), \nu) \bigr]_A
		&=
		\mleft( 
			- K \mleft( H(\mu, \nu) \mright),
			H(\mu, \nu)
		\mright) ~ ,
		\end{align*}
		\textit{where}
		\begin{equation*}
		H(\mu, \nu)
		=
		t_{\basicc}(\mu, \nu)
			+ \zeta \mleft( K(\mu), K(\nu) \mright)
		\end{equation*}
		\textit{for all $\mu, \nu \in \Gamma(E)$.}
		\end{theorems}
		
		Observe the similarity with matched pairs of Lie algebroids, \cite{laurent2008holomorphic}, for which $\zeta$ is zero and thus making use of an $F$-representation on $E$, while here only $\basicc$ is flat; but we on the other hand assume the existence of $K$, {\color{Cred} and $\basicc$ is defined by $\nabla$ and in general not independent of it}. 
		
		\cite{Fischer:2024vak} is written in the BRST formalism so that it might be difficult for the unfamiliar reader to check the reference; thus, {\color{Cred}we will provide an alternative proof, but first} a short explanation {how to use the reference}: Lie algebroids are equivalent to differential graded vector bundles of degree 1 equipped with a cohomological vector field, in particular the anchor and the Lie bracket can be read of the differential as in \cite[interpreting the Weil differential of Equation (2.55a) as a cohomological vector field; components along the base manifold encodes the anchor, components along the fibres the Lie bracket]{Fischer:2024vak}.
		
		\begin{proof}[Proof of \ref{thm:CanonEx}]
		{\color{Cred}We will prove this by using} Mackenzie's studies as {\color{Cred}reiterated} in \ref{sec:basics}, and {\color{Cred}by using} \ref{thm:ActionAlgebroidStructure}. In order to proceed in this manner, let us first rewrite the structure of $A$ a bit. Observe that the graph ${{\rm Graph}}(-K) \coloneqq \{(-K(\mu), \mu) ~ | ~ \mu \in E\}$ of $-K$ is a BLA by \ref{thm:ActionAlgebroidStructure}; {\color{Cred}it is straightforward to check that the proposed Lie bracket of $A$ restricts to $H$ by restricting to ${{\rm Graph}}(-K)$}. Moreover, it is the kernel of $\scD \colon A \to F$, $(X, \mu) \mapsto X + K(\mu)$, which is clearly surjective. {\color{Cred}For covenience,\footnote{So that we can write ``morphism of Lie algebroids'' instead of ``morphism of almost Lie algebroids'' and so on.} let us also assume for a moment that $A$ is indeed a Lie algebroid, then $\scD$} is also a morphism of Lie algebroids:
		\begin{equation*}
		\rho_F \circ \scD
		=
		\rho_A ~ ,
		\end{equation*}
		and
		\begin{equation*}
		\scD\mleft( \mleft[ (X, \mu) ~ , (Y, \nu) \mright]_A \mright)
		=
		\mleft[ \scD(X, \mu), \scD(Y, \nu) \mright]_F ~
		\end{equation*}
		for all $(X, \mu), (Y, \nu) \in \Gamma(A)$, where we made use of $K$ being a morphism of Lie algebroids. 
		Henceforth, we have the following short exact sequence of Lie algebroids
		\begin{equation}\label{eq:TypicalShortExactSequence}
			\begin{tikzcd}
				{\rm Graph}(-K) \arrow[r, hook]& A \arrow[r,"\scD", two heads] & F ~ ,
			\end{tikzcd}
		\end{equation}
		where the embedding of ${\rm Graph}(-K)$ into $A$ is the canonical one. It admits a {\color{Cred}strict YM $F$-connection $\hat{\nabla}$} coming from a splitting of $\scD$; let us choose the canonical splitting {\color{Cred} $\chi \colon F \to A$, $X \mapsto (X, 0)$, then $\scD \circ \chi = {\rm id}_F$ and the induced connection} $\hat{\nabla}$ is given as {\color{Cred}usual by}
		\begin{align*}
		\hat{\nabla}_X\mleft( -K(\mu), \mu \mright)
		&=
		\mleft[ (X, 0) ~ , \mleft(-K(\mu), \mu\mright) \mright]_A
		\\
		&=
		\mleft(  
			- K \mleft( \nabla_X \mu - \zeta\mleft( X, K(\mu) \mright) \mright) ~,
			\nabla_X \mu - \zeta\mleft( X, K(\mu) \mright) 
		\mright)
		\\
		&=
		\mleft(  
			- K \mleft( \nabla^\zeta_X \mu \mright) ~,
			\nabla^\zeta_X \mu 
		\mright)
		\end{align*}
		for all $X \in \Gamma(F)$ and $\mu \in \Gamma(E)$.
		That is, $\hat\nabla \circ \iota = \iota\circ \nabla^\zeta$, where $\iota\colon E \to \mathrm{Graph}(-K)$, $\iota(\nu) \coloneqq (-K(\nu), \nu)$. Observe that there is a more convenient way to write the short exact sequence \eqref{eq:TypicalShortExactSequence}: $\iota$ satisfies
		\begin{equation*}
		\mleft[ \iota(\mu) ~ , \iota(\nu) \mright]_{{\rm Graph}(-K)}
		=
		\iota\mleft( H(\mu, \nu) \mright) ~ ,
		\end{equation*}
		where we denote the bracket on ${\rm Graph}(-K)$ with the corresponding subscript. It follows that $\iota$ is an isomorphism of BLAs $E_H \cong {\rm Graph}(-K)$ (as it is already a canonical isomorphism of vector bundles), which implies that the short exact sequence \eqref{eq:TypicalShortExactSequence} can be also written as
		\begin{equation}\label{eq:TheShortExactSequence}
			\begin{tikzcd}
				E_H \arrow[r, hook, "\iota"]& A \arrow[r,"\scD", two heads] & F ~ .
			\end{tikzcd}
		\end{equation}
		Now {\color{Cred}we are finally able to} rewrite the structure on $A$, that is,
		\begin{align}\label{eq:BracketsOfBothAStructures}
		\mleft[ (X, \mu) ~ , (Y, \nu) \mright]_A
		&=
		\mleft[ (\scD(X, \mu) ~, 0) + \iota(\mu), (\scD(Y, \nu) ~, 0) + \iota(\nu) \mright]_A
		\nonumber\\
		&=
		\mleft( \mleft[ \scD(X, \mu), \scD(Y, \nu) \mright]_F ~ , 0 \mright)
		\nonumber\\
		&\hspace{.6cm}
			+ \mleft[ \iota(\mu), \iota(\nu) \mright]_A
			+ \hat\nabla_{\scD(X, \mu)} \iota(\nu)
			- \hat\nabla_{\scD(Y, \nu)} \iota(\mu)
			+ \iota\mleft(\zeta\mleft( \scD(X, \mu), \scD(Y, \nu) \mright) \mright)
		\nonumber\\
		&=
		\mleft( \mleft[ \scD(X, \mu), \scD(Y, \nu) \mright]_F ~ , 0 \mright)
		\nonumber\\
		&\hspace{1cm}
			+ \iota\mleft( 
				H(\mu, \nu) 
				+ \nabla^\zeta_{\scD(X, \mu)} \nu
				- \nabla^\zeta_{\scD(Y, \nu)} \mu
				+ \zeta\mleft( \scD(X, \mu), \scD(Y, \nu) \mright)
			\mright) ~ .
		\end{align}
		In particular, by making use of $\scD(X - K(\mu), \mu) = X$,
		\begin{align*}
		\mleft[ (X, 0) + \iota(\mu) ~ , (Y, 0) + \iota(\nu) \mright]_A
		&=
		\mleft[ (X - K(\mu), \mu ) ~ , (Y - K(\nu), \nu) \mright]_A
		\\
		&=
		\mleft( \mleft[ X, Y \mright]_F ~ , 0 \mright)
			+ \iota\mleft( 
				H(\mu, \nu) 
				+ \nabla^\zeta_{X} \nu
				- \nabla^\zeta_{Y} \mu
				+ \zeta\mleft( X, Y \mright)
			\mright) ~ .
		\end{align*}
		
		This is precisely the construction as in \ref{sec:basics}, based on the short exact sequence \eqref{eq:TheShortExactSequence}, and {\color{Cred}this proves that $A$ is indeed a Lie algebroid. Spelling things out: Observe that the proposed Lie bracket on $A$ is bilinear, antisymmetric, and satisfies the Leibniz rule w.r.t.\ $\rho_A$ by construction; furthermore, $\rho_A$ is a homomorphism of brackets (making use of $\rho_F \circ K = \rho_E$), so that $A$ is an \textit{almost} Lie algebroid. The previous identities hold regardless of whether or not $\mleft[ \cdot, \cdot \mright]_A$ satisfies the Jacobi identity, in particular the previous discussion holds in the category of almost Lie algebroids, and the last equation can be shown. The Jacobi identity now follows by the fact that $(\nabla^\zeta, \zeta)$ is a strict YM $F$-connection (\ref{thm:ActionAlgebroidStructure}), in the same fashion as in \cite[Theorem 7.3.7]{0521499283}.}
		\end{proof}
		
		{\color{Cred} Essentially, we observed that $A = F\bowtie_{(\nabla, \zeta)} E$ sits in a short exact sequence of Lie algebroids \eqref{eq:TheShortExactSequence}, and thus the Lie algebroid structure of $A$ is isomorphic to the structure provided in the discussion around \eqref{eq:MackAlgebroid}, so that we can on one hand concisely summarise the Lie algebroid structure on $A$ and on the other hand complement \ref{thm:CanonEx} with \ref{thm:ActionAlgebroidStructure} in the following corollary:}
		
		\begin{corollaries}{Comparing the two structures on $A$}{ComparingTwoIsomOnA}
		\color{Cred} Given the notation of \ref{thm:CanonEx} and its proof: We have an isomorphism of algebroids given by 
		\begin{align*}
		F \bowtie_{(\nabla, \zeta)} E &\to F\ltimes_{(\nabla^{\zeta}, \zeta)} E_{H}~,\\
		(X, \mu) &\mapsto \mleft( \scD(X, \mu), \mu \mright)~,
		\end{align*}
		where $\scD \colon F \bowtie_\nabla E \to F$ is the Lie algebroid morphism given as $\scD(X, \mu) \coloneqq X + K(\mu)$. In particular,
		\begin{equation*}
		\mleft[ \mu, \nu \mright]_E
			+ \nabla_X \nu - \nabla_Y \mu
			+ \zeta(X, Y)
		=
		H(\mu, \nu)
			+ \nabla^{\zeta}_{\scD(X, \mu)}\nu
			- \nabla^{\zeta}_{\scD(Y, \nu)}\mu
			+ \zeta\mleft( \scD(X,\mu), \scD(Y, \nu) \mright)
		\end{equation*}
		for all $(X, \mu), (Y, \nu) \in \Gamma(F \oplus E)$.
		\end{corollaries}
		
		{\color{Cred} 
		
		\begin{proof}
		The identity for the Lie brackets is just a consequence of comparing the second components in Eq.\ \eqref{eq:BracketsOfBothAStructures}. Using this, or by applying the proposed isomorphism of Lie algebroids on both sides of Eq.\ \eqref{eq:BracketsOfBothAStructures}, now shows that the proposed isomorphism is a morphism of Lie brackets, while all other needed properties follow immediately by definition.
		\end{proof}
		
		Let us keep the notation of the previous two statements and their proofs; in the beginning of the proof we assumed first that $A$ is indeed a Lie algebroid, and we will now extend this discussion:} Given the splitting $A = F \oplus E$, we have a canonical short exact sequence of vector (!) bundles
		\begin{equation}\label{eq:LieAlgVBsequence}
			\begin{tikzcd}
				E \arrow[r, "\hat\iota", hook]& A \arrow[r, "\psi", two heads] & F ~ ,
			\end{tikzcd}
		\end{equation}
		The canonical projection $\psi$ onto $F$ is in general not even a morphism of anchored vector bundles, because this would otherwise imply that $E$ is a bundle of Lie algebras as explained at the beginning of \ref{sec:basics}; that behaviour of $\psi$ can be easily confirmed by using the structure provided in \ref{thm:CanonEx}. However, the canonical embedding $\hat\iota \colon E \to A$, $\mu \mapsto (0, \mu)$, is in fact a morphism of Lie algebroids due to the fact that
		\begin{equation*}
		\rho_A \circ \hat \iota = \rho_E ~ ,
		\end{equation*}
		and
		\begin{equation*}
		\mleft[ \hat\iota(\mu), \hat\iota(\nu) \mright]
		=
		\hat\iota\mleft( \mleft[ \mu, \nu \mright]_E \mright) ~.
		\end{equation*}
		Moreover, as previously, let us again take the canonical lift $\chi \colon F \to A$, $X \mapsto (X, 0)$, but now acting on $E$ via the short exact sequence \eqref{eq:LieAlgVBsequence}. That is,
		\begin{equation*}
		\mleft[ \chi(X), \hat\iota(\mu) \mright]_A
		=
		\mleft( 
			-\basicc_\mu X,
			\nabla_X \mu
		\mright) ~.
		\end{equation*}
		By applying $\scD$ on both sides one recovers the definition of the basic connection on $F$ (making use of the fact that $\scD$ is a morphism of Lie algebroids). One also has
		\begin{align*}
		\nabla_X \mu
		&=
		\hat\chi \mleft( \mleft[ \chi(X), \hat\iota(\mu) \mright]_A \mright) ~ ,
		\\
		\basicc_\mu X
		&=
		\psi\mleft( \mleft[ \hat\iota(\mu), \chi(X) \mright]_A \mright) ~ ,
		\end{align*}
		where $\hat\chi\colon A \to E$ is the canonical projection on $E$, in fact it is naturally the retro-splitting of $\chi$; recall that those are uniquely defined by $\hat\chi \circ \hat\iota = {\rm id}_E$ and 
		\begin{equation*}
		{\rm id}_A = \hat\iota \circ \hat\chi + \chi \circ \psi ~.
		\end{equation*}
		Since $\psi$ is not a morphism of anchored vector bundles, one achieves $\basicc$ as a nontrivial $E$-connection on $F$, additionally to the typical construction of $\nabla$. {\color{Cred}Actually, we can write $\psi = \scD - K \circ \hat\chi$ as a perturbance of $\scD$, which is also linked to the minimal coupling in physics}. If $\psi$ is a morphism of anchored vector bundles, then $\basicc$ is trivial as one expects; however, $\basicc$ is certainly flat by assumption, while we have for the curvature $R_\chi$ of $\chi$ that
		\begin{equation*}
		R_\chi(X,Y) 
		= 
		\mleft[ \chi(X), \chi(Y) \mright]_A
			- \chi \mleft( \mleft[ X, Y \mright]_F \mright)
		=
		\iota\mleft( \zeta(X, Y) \mright) ~ .
		\end{equation*}
		In particular, while $R_\chi$ is not in the kernel of $\psi$ in general because $\psi$ is not a morphism of algebroids, we have 
		\begin{equation*}
		\scD \circ R_\chi = 0 ~ ,
		\end{equation*}
		which we already knew because $\chi$ is also a splitting of \eqref{eq:TheShortExactSequence} which \textit{is} a short exact sequence of Lie algebroids. Therefore, if $\nabla^\zeta$ is flat, then it describes a Lie algebroid morphism embedding $F$ into $A$; as usual, the existence of a flat splitting of \eqref{eq:TheShortExactSequence} is locally given.
		
		In total we have the following commuting diagram:
		\begin{equation}\label{eq:TotalDiag}
			\begin{tikzcd}[column sep=small]
			E \arrow[dd, "K", swap] \arrow[rd, "\hat\iota", hook] &  & E_H \arrow[ld, "\iota", swap, hook] \arrow[dd, "-K", dotted]
			\\
			& A \arrow[rd, "\psi", swap, two heads, dotted] \arrow[dl, "\scD", two heads] & 
			\\
			F & & F
			\end{tikzcd}
		\end{equation}
		where the dotted lines are only vector bundle morphisms in general. This \textbf{sandglass sequence coupling $E$ and $F$ {\color{Cred}following $K$}} will provide the starting point of the second paper where we will discuss the obstruction behind the Lie algebroid structure on $A$ making use of strict covariant $K$-adjustments, and where we will clarify the notion of splittings and their changes {\color{Cred}in order to generalise BLA ii) of the introduction and possible related obstructions}; all of this being answered with the tools of curved Yang-Mills-Higgs theories.
		
		\section{Examples}\label{sec:Examples}
		
		Let us conclude this paper with several examples of sandglass sequences; however, most examples with strict adjustments and Yang-Mills connections were already presented in \cite{Fischer:2020lri, Fischer:2021glc, Fischer:2022sus, Fischer:2024vak}, so that we will not repeat {\color{Cred} and discuss} those examples in detail, see those references for elaborated details instead. {\color{Cred} Let us instead give a short introduction of those examples before we highlight certain types of examples, including new ones.
		}

		Here we will only introduce the abstract idea of the most important examples, and afterwards we will turn to a new class of examples motivated by \cite{Fischer:2401.05966}, concluding with crossed modules and their obstruction. But let us start with the obvious, {\color{Cred}keeping the notation of the previous section}:

		\subsection{The edge cases}
		
		\begin{examples}{Atiyah sequences}{BLAsasExample}
		Of course, we recover all the typical Atiyah sequences, that is, $E$ itself being an BLA: Assume $K = 0$, then {\color{Cred}$E$ is a BLA}, $E = E_H$ by \ref{thm:ActionAlgebroidStructure}, {\color{Cred} and strict $0$-adjustments are strict YM $F$-connections (\ref{cor:KEqual0IsYM}). Those are precisely those coming from a splitting of a short exact sequence of Lie algebroids with $E$ as a kernel and $F$ as image:}
		
		We also have $\hat \iota = \iota$ and $\scD(X, \mu) = \psi (X, \mu) + (K \circ \hat \chi)(X, \mu)$, that is, $\scD = \psi$ for a vanishing $K$. Then also $\nabla = \nabla^\zeta$ is a strict multiplicative Yang-Mills $F$-connection, $\basicc \equiv 0$ (on $F$) and $\basicc = [\cdot, \cdot]_E$ (on $E$), such that we recover the construction of \ref{sec:basics}. 
		\end{examples}
		
		{\color{Cred}Unit octonions are in fact an example, too:}
		
		\begin{examples}{Tangent algebroids}{TangentAlgoids}
		\tcbsubtitle{\cite{Fischer:2021glc,Fischer:2024vak}}
		\color{Cred} Let us now assume that $E = F = \rmT M$ and $K = \mathrm{id}_{\rmT M}$. Cartan connections $\nabla$ are in this case 1:1 to the flat connection $\basicc$ by definition, in particular $\mleft( \basicc \mright)^{{\rm bas}} = \nabla$. Recall \ref{rem:PossibleChoiceOfZetaKK}, a primitive \textit{always} exists in this case, in particular, $\zeta = t_\nabla$ is a valid choice which makes $(\nabla, \zeta)$ a strict adjustment; see \cite{Fischer:2024vak}. Choosing such a primitive induces an abelian strict LAB\footnote{\color{Cred}As discussed in \ref{rem:LABalongOrbits}, the transitivity of $F = \rmT M$ implies that the strict BLA structure is actually an LAB.} $\rmT M_H$, however, one can vary $\zeta$ with forms preserving the compatibilities to induce other LAB structures.
		
		As a toy example: $M = \mathbb{S}^7$, the unit octonions, only admits curved Cartan connections whose curvature is directly linked to the lack of associativity in the octonionic algebra. In general every parallelizable, compact, connected and simply connected smooth manifold $M$ can be equipped with a strict adjustment for which the adjustment is flat if and only if $M$ admits a Lie group structure; see also \ref{thm:ObstrCrossedModules} later.
		\end{examples}
		
		\subsection{Action Lie algebroids}
		
		The next examples will be based on the pullback of Lie algebroid connections and action algebroids; let us first quickly recapitulate the previously mentioned relationship to action algebroids, recall \ref{thm:ActionAlgebroidStructure}, {\color{Cred}\ref{ex:RecoverTrivialActionAlgoid}, and in particular the flat example as given in \ref{cor:RevisitObstrActAlgoid}}:
		
		\begin{examples}{Action Lie algebroid, the trivial and flat example}{ClassActionAlg}
		\color{Cred}Let $\frg$ be a Lie algebra acting on a manifold $M$. Then we define $E \coloneqq M \times \frg$ to be the canonical action algebroid, $F \coloneqq \rmT M$, $K \coloneqq \rho_E$, while we choose the canonical flat connection $\nabla^0$ as the strict adjustment on $E$ where the primitive $\zeta$ is chosen to be zero, and thus \color{Cred}$\nabla^\zeta = \nabla^0$; moreover, w.r.t.\ constant sections $\mu, \nu$ of $E$,
		\begin{equation*}
		H(\mu, \nu)
		=
		t_{\nabla^{0,\mathrm{bas}}}(\mu, \nu)
			+ \zeta \mleft( K(\mu), K(\nu) \mright)
		=
		- t_{\nabla^0_{\rho_E}}(\mu, \nu)
		=
		\mleft[\mu, \nu\mright]_E
		=
		\mleft[\mu, \nu\mright]_{\frg}~,
		\end{equation*}
		by making use of the action algebroid structure and $t_{\basicc} = -t_{\nabla_K}$. That is, as expected, the strict LAB $E_H$ is the canonical LAB structure on $M \times \frg$. By \ref{thm:ActionAlgebroidStructure} we also have
		\begin{align*}
		\mleft[ \mu, \nu \mright]_E
		&=
		H(\mu, \nu)
			+ \nabla^\zeta_{K(\mu)} \nu
			- \nabla^\zeta_{K(\nu)} \mu
			+ \zeta \mleft( K(\mu), K(\nu) \mright)
		\\
		&=
		\mleft[ \mu, \nu \mright]_\frg
			+ \nabla^0_{\rho_E(\mu)} \nu
			- \nabla^0_{\rho_E(\nu)} \mu
		\end{align*}
		for all $\mu, \nu \in \Gamma(E)$, which is just the typical form of the action algebroid bracket. 
		
		Let us now turn to $A = \rmT M \bowtie_{(\nabla^0,\zeta \equiv 0)} E$ as given in \ref{thm:CanonEx} for which we rewrite the bracket a bit by using $\rho_A(X, \mu) = \scD(X, \mu) = X + \rho_E(\mu)$ for $(X, \mu) \in A = F \oplus E$:
		\begin{align*}
		\mleft[ (X, \mu), (Y, \nu) \mright]_A
		&\coloneqq
		\Bigl( 
			\mleft[ \rho_A(X, \mu), \rho_A(Y, \nu) \mright]_{\rmT M}
			- \rho_E \bigl( 
				\mleft[ \mu, \nu \mright]_\frg + \nabla_{\rho_A(X, \mu)} \nu - \nabla_{\rho_A(Y, \nu)} \mu
			\bigr)~,
		\\
		&\hspace{1cm}
			\mleft[ \mu, \nu \mright]_\frg + \nabla_{\rho_A(X, \mu)} \nu - \nabla_{\rho_A(Y, \nu)} \mu
		\Bigr)
		\end{align*}
		for all $(X, \mu), (Y, \nu) \in \Gamma(A)$; alternatively one can write the first component using $\nabla^{0, \mathrm{bas}}$ as in \ref{thm:CanonEx}, or one makes use of \ref{cor:ComparingTwoIsomOnA}.
		\end{examples}
		
		Now a short reminder of how pullbacks of Lie algebroid connections work as for example proven in \cite[Corollary 3.5.7]{Fischer:2021glc}; terms starting with $\phi^*$ denote the pullback of vector bundles and associated pullback constructions.
		
		\begin{corollaries}{Pullbacks of Lie algebroid connections}{PullbackofFConn}
		Let $F_i \to M_i$ ($i \in\{1,2\}$) be two Lie algebroids over smooth manifolds $M_i$, $E \to M_2$ a vector bundle, and ${}^{F_2}\nabla$ an $F_2$-connection on $E$. Also fix an anchor-preserving vector bundle morphism $\xi\colon F_1 \to F_2$ over a smooth map $\phi\colon M_1 \to M_2$. Then there is a unique $F_1$-connection $\phi^*\mleft( {}^{F_2}\nabla \mright)$ on $\phi^*E$ with
		\begin{equation*}
		\mleft(\phi^*\mleft( {}^{F_2}\nabla \mright)\mright)_Y (\phi^*\mu)
		=
		\phi^*\mleft(
			{}^{F_2}\nabla_{\xi(Y)} \mu
		\mright)
		\end{equation*}
		for all $\mu \in \Gamma(E)$ and $Y \in \Gamma(F_1)$.
		\end{corollaries}
		
		As pointed out in \cite[\S 7.3]{Fischer:2024vak}, given an algebroid $E \to M$, equipped with a strict covariant $K$-adjustment $\nabla$, it is possible to make a certain pullback of $\nabla$: Given an $E$-action {\color{Cred}$\Gamma(E) \to \frX(N)$, $\mu \to \overline\mu$,} on a smooth manifold $N$ along a submersion $\phi \colon N \to M$ as a moment map, one can take the pullback of $F$ as Lie algebroid, denoted by $\phi^!F$, coming with a canonical morphism of Lie algebroids $\xi\colon \phi^!F \to F$ over $\phi$, {\color{Cred} which is surjective if viewed as a base-preserving vector bundle morphism $\Phi^!F \to \Phi^*F$}; this is due to the fact that $\phi$ is transverse to the anchor of $F$ such that we can apply the standard construction for pullbacks of Lie algebroids, that is, $\phi^!F$ consists of elements in $(Y, \eta) \in \phi^*F \oplus \rmT N$ such that $\mleft(\phi^*\rho_F\mright)(Y) = \rmD \phi(\eta)$.
		
		Regarding $E$ however, we only take its pullback as vector bundle $\phi^*E$ and equip it with the canonical action Lie algebroid structure, {\color{Cred}that is,\footnote{\color{Cred}In the following we abuse the notation of what we mean with $\mu$ a bit, for simplified notation. It should be clear by context what type of element it is.} following \cite[Proposition 4.1.2]{0521499283}, $\rho_{\Phi^*E}(p, \mu) = \overline{\mu}_p$ for all $(p, \mu) \in \Phi^*E$, and the Lie bracket on its sections is uniquely given by
		\begin{equation*}
		\mleft[ \Phi^*\mu, \Phi^*\nu \mright]_{\Phi^*E}
		=
		\Phi^*\mleft(\mleft[ \mu, \nu \mright]_{E}\mright)
		\end{equation*}
		for all $\mu, \nu \in \Gamma(E)$.} The morphism $K$ naturally extends to a Lie algebroid morphism $\phi^!K \colon \phi^*E \to \phi^!F$, $\mu \mapsto \mleft( \mleft(\phi^*K\mright)(\mu), \rho_{\phi^*E}(\mu) \mright)$, {\color{Cred} which is well-defined due to} the fact that the definition of Lie algebroid actions implies 
		\begin{equation*}
		\mleft(\phi^*\rho_F \circ \phi^*K \mright)\mleft( \mu \mright) 
		=
		\mleft(\phi^*\rho_E\mright)(\mu)
		= 
		\rmD \phi\mleft(\rho_{\phi^*E}(\mu)\mright) ~ .
		\end{equation*}
		Then the action algebroid $\phi^*E$ comes with a natural strict $\phi^!K$-adjustment, that is, $(\phi^*\nabla, \zeta^\prime)$ is a strict $\phi^!K$-adjustment, where $\phi^*\nabla$ is the pullback of $\nabla$ as in \ref{cor:PullbackofFConn} along $\xi$, that is, $\phi^*\nabla$ is the canonical $\phi^!F$-connection on $\phi^*E$, and where $\zeta^\prime$ is uniquely given by
		\begin{equation*}
		\zeta^\prime\mleft(\phi^!X, \phi^!Y\mright)
		=
		\phi^*\mleft(\zeta\mleft( X, Y \mright)\mright)
		\end{equation*}
		for all $X, Y \in \Gamma(F)$; here $\phi^!X, \phi^!Y \in \Gamma(\phi^!F)$ are any sections which project to $X, Y$ under $\xi$, respectively; for this it is essential to observe that setions of $\phi^!F$ are generated precisely by sections of the form $\phi^!X$ because $\xi$ is a surjective morphism {\color{Cred}$\Phi^! F \to \Phi^*F$}. The mentioned reference shows the above for $F$ being a tangent algebroid, but the proof is precisely the same for arbitrary $F$ by observing that we have
		\begin{align*}
		\mleft(\xi \circ \phi^!K\mright)(\phi^*\mu)
		&=
		\phi^*\mleft( K(\mu) \mright)
		\\
		\mleft( \phi^*\nabla \mright)_{\phi^!X} \phi^*\mu
		&=
		\phi^*\mleft( \nabla_X \mu \mright)
		\end{align*}
		for all $\mu \in \Gamma(E)$ and $X \in \Gamma(F)$, in particular also $\mleft(\phi^!K\mright)(\phi^*\mu) = \phi^!\mleft( K(\mu) \mright)$. The proof that $(\phi^*\nabla, \zeta^\prime)$ is a strict covariant $\phi^!K$-adjustment is then precisely the same as in \cite[\S 7.3]{Fischer:2024vak} since it uses typical pullback arguments.
		
		\begin{examples}{Pullbacks of sandglass sequences along actions}{Pullback}
		In particular, if there is a sandglass sequence coupling $E \to M$ and $F \to M$ {\color{Cred} via a strict $K$-adjustment $(\nabla,\zeta)$}, then there is also one coupling $\phi^*E \to N$ and $\phi^!F \to N$ {\color{Cred} via $(\Phi^*\nabla,\zeta')$}. If $F = \rmT M$ and $K = \rho_E$, then $\phi^! F \cong \rmT N$ naturally by projection and thus $\phi^!K \cong \rho_{\phi^*E}$.
		
		{\color{Cred}Due to the fact that all of the involved constructions are usual pullback arguments it is straight-forward to check that the strict BLA of $\Phi^*E$ is the pullback-BLA $\Phi^*(E_H)$ with field of Lie brackets $\Phi^*H$, and thus we have the canonical sandglass sequence
		\begin{equation*}
			\begin{tikzcd}[column sep=small]
			\Phi^*E \arrow[dd, "\Phi^!K", swap] \arrow[rd, hook] &  & \Phi^*(E_H) \arrow[ld, swap, hook] \arrow[dd, "-\Phi^!K", dotted]
			\\
			& \Phi^!F \bowtie_{(\Phi^*\nabla,\zeta')} \Phi^*E \arrow[rd, swap, two heads, dotted] \arrow[dl, two heads] & 
			\\
			\Phi^!F & & \Phi^!F
			\end{tikzcd}
		\end{equation*}
		Last, we have a canonical isomorphism of Lie algebroids 
		\begin{align*}
		\Phi^!F \bowtie_{(\Phi^*\nabla,\zeta')} \Phi^*E 
		&\to 
		\Phi^!\mleft( F\bowtie_{(\nabla, \zeta)} E \mright)~,
		\\
		(p,Y,\eta, \mu)
		&\mapsto
		\mleft(p, Y, \mu, \eta + \rho_{\Phi^*E}(p, \mu)\mright)~,
		\end{align*}
		where $p \in N$, $Y\in F_{\Phi(p)}$, $\eta \in \rmT_pN$, and $\mu \in E_{\Phi(p)}$; this also follows just by construction and by definition of $\Phi^!K$. This allows to express the unlabelled arrows of this pullback sandglass sequence as pullback maps of \eqref{eq:TotalDiag}, similar to the construction of $\Phi^!K$. However, since all of this is straight-forward to show, but a bit tedious, and since we will not need the pullback that explicit, we leave the explicit proof as an exercise to the reader.}
		\end{examples}
		
		\begin{examples}{Transitive algebroid acting on normal bundle}{transitiveE}
		Of a particular interest might be a transitive algebroid $E$ acting on a normal bundle of $M$ in $N$ preserving the 0-section ($=M$), where we set $F = \rmT M$ and $K = \rho_E$; $M$ is the a leaf of the singular foliation in the normal bundle generated by the anchor of $\phi^*E$. As pointed out in \cite{Fischer:2024vak}, such $E$ locally admit strict adjustments, implying the local existence of the mentioned sandglass sequences; the reference further assumed faithfulness of the action, however, this was done for other reasons and not needed for the local existence which just works for any action, as long as $\phi$ is a submersion.
		\end{examples}
		
		From a practical point of view it is often much easier to study BLAs $E$ (with $K=0$) in order to produce an adjustment for the action algebroid structure on $\phi^*E$, {\color{Cred} because there are many known examples of short exact sequences of Lie algebroids and it allows to produce adjustments on Lie algebroids which are not a BLA}. 
		
		\begin{examples}{Atiyah sequences producing sandglasses along $\Phi$}{AtiyahAction}
		\color{Cred} Given an ordinary short exact sequence of Lie algebroids with $E$ as kernel, naturally coming with strict 0-adjustments (recall \ref{ex:BLAsasExample}), assume that the BLA $E$ acts on a manifold $N$ along a submersion $\Phi \colon N \to M$. Then it produces a sandglass sequence as in \ref{ex:Pullback} for which $\Phi^!K$ is in general non-trivial, and $\Phi^*E$ is a Lie algebroid which is in general not a BLA.
		
		As a cooking recipe: Any principal bundle $P \to M$ with structural Lie group $G$ comes with its canonical Atiyah sequence which induces a strict 0-adjustment on the adjoint bundle ${\rm ad}(P) \coloneqq (P \times \frg)/G$, where $G$ acts on its Lie algebra $\frg$ by the adjoint representation; the canonical strict 0-adjustments available on ${\rm ad}(P)$ are precisely the adjoint connections induced by an Ehresmann connection on $P$. Via an action of ${\rm ad}(P)$ along $\Phi$ one achieves an action algebroid structure on $\Phi^*{\rm ad}(P)$ with the canonical strict $\Phi^!K$-adjustment by pulling back the adjoint connection, where $\mleft(\Phi^!K\mright)(p,\mu) = (0, \overline\mu_p)$. For example, the adjoint bundle of the Hopf fibration $P \coloneqq \mathbb{S}^7 \to \mathbb{S}^4$ canonically acts on the associated bundle $(P \times \mathbb C^2)/\mathbb S^3$, leading to an action Lie algebroid which only admits curved strict adjustments; see \cite[Example 7.25]{Fischer:2024vak}.
		\end{examples}
		
		As discussed above, the strict LAB of the action Lie algebroid $\phi^*E$ in \ref{ex:AtiyahAction} is canonically given by the pullback BLA structure coming from $E$ {\color{Cred}due to $E = E_H$} which highlights the duality of action algebroids and BLA structures on the same bundle as in \ref{thm:ActionAlgebroidStructure}. {\color{Cred}This allows to easily construct Lie algebroid examples by starting with BLAs; see \cite{Fischer:2022sus, Fischer:2401.05966, Fischer:2024vak} for partial classifications and discussions regarding examples not admitting flat adjustments. In particular, this approach allows to easily construct strict $\rho_E$-adjustments by looking at ``ordinary'' Atiyah sequences ($F = \rmT M$), and then adding an action.} Thus, let us {\color{Cred}now highlight} interesting LABs $E$ {\color{Cred}which one can use to construct such strict adjustments}.
		
		Following \cite{Fischer:2401.05966}, there are LABs found to be in anchored vector bundles, Poisson geometry and so on, {in general, whenever one has a ``splitting theorem'', inducing a unique transverse structure}:
		
		We now follow \cite{bursztyn2019splitting}. Given an anchored vector bundle $W \to M$, its anchor induces a singular foliation on $M$. Fix a leaf $L$ of this foliation and consider submanifolds in $M$ transverse to the foliation and intersecting $L$ trivially, then it is a well-known fact that $W$ restricted to those transverse submanifolds is again an anchored vector bundle with imprinted singular foliation; these restrictions of $W$ are also called \textit{transverse structures of $W$ along $L$}. Locally, the transverse structure at two different points in $L$ are isomorphic to each other and this is due to the fact that there is a horizontal lift of vectors in $L$ to a vector preserving the transverse structure. In \cite{Fischer:2401.05966} we extended this idea and showed that this leads to a strict multiplicative Yang-Mills $\rmT L$-connection:
		
		\begin{examples}{Locally split structures}{LocallSplitLABs}
		Either assume a formal setting, or assume that $L$ is an embedded leaf and that $M$ is the normal bundle of $L$. Each fibre of the normal bundle is imprinted with a transverse structure of $W$, all isomorphic to each other such that we denote the structural transverse structure by $\tau(W_*)$, where $*$ denotes any fixed point on $L$. 
		
		The group of automorphisms $\rm{Aut}(\tau(W_*))$ of $\tau(W_*)$ naturally forms a Lie group bundle over $L$, with induced LAB $E$. Here automorphism means that it comes with extra structure, depending on the structure of $W$ as in \cite{bursztyn2019splitting}. If $W$ is an anchored vector bundle, then we mean automorphisms of anchored vector bundles; if $W$ is a Lie algebroid, then we mean automorphisms of Lie algebroids; if $W$ is a Poisson Lie algebroid, then we mean Poisson automorphisms, and so on.
		
		By \cite{bursztyn2019splitting} there is a natural parallel transport of the transverse structures along $L$ which is an isomorphism of transverse structures, in particular over closed loops it has values in $\rm{Aut}(\tau(W_*))$. Naturally, the curvature thence has values in the connected component of $\rm{Aut}(\tau(W_*))$ around the identity, and as in \cite{Fischer:2401.05966} this is equivalent to the notion of multiplicative Yang-Mills $\rmT L$-connections on the group bundle and so also on $E$, {\color{Cred}infinitesimally}. Strictness naturally comes by the fact that the group bundle naturally sits in a short exact sequence as the isotropy bundle of the groupoid of automorphisms of transverse structures between different points.
		
		As a special case: If $W$ is a Poisson Lie algebroid, then the horizontal lift of this strict multiplicative Yang-Mills $\rmT L$-connection lifts to Poisson vector fields, while the curvature has values in Hamilton vector fields. That is, the parallel transport has values in the group of Poisson isomorphisms of the transverse structures, and the curvature has values in the group of Hamiltonian diffeomorphisms.
    \end{examples}
		
		\begin{remark}
		Besides the mentioned references, \cite{meinrenken2021integration} is also a very useful reference for the reader unfamiliar with such constructions. The second paper will make the construction {\color{Cred} above} even clearer once we introduced the obstruction; but we can already observe that the parallel transport followed by the quotient map of $\rm{Aut}(\tau(W_*))$ over its identity component is flat. This is what Mackenzie called coupling and what is called \textit{outer holonomy} in \cite{Fischer:2401.05966}. We will generalise this and its relation to sequences, making use of the fact that the curvature of $\nabla$ satisfies an exactness condition.
		\end{remark}
		
		{\color{Cred}
		\subsection{Crossed modules}
		
		Even though this series of papers is motivated by generalising $E$ to a Lie algebroid, assuming that $E$ is an LAB does not mean we just have ordinary definitions: $K$ can still be non-trivial. An important example are actually crossed modules which naturally fit in the setting of sandglass sequences; we start with the definition of crossed modules of Lie algebroids as introduced in \cite{androulidakis2005crossed} (and also following \cite{laurent2008obstruction}), while we drop the transitivity of $F$ as in the references so that it is obvious how to recover the ``ordinary'' definition of crossed modules of Lie algebras/(trivial) LABs.
		
		\begin{definitions}{Crossed modules of Lie algebroids}{CrossedAlgoids}
		\tcbsubtitle[before skip=\baselineskip]{\cite[Definition 2.2]{androulidakis2005crossed}}
		{\color{Cred} Let $E \to M$ be an LAB over a smooth manifold $M$, $F \to M$ a Lie algebroid, $K \colon E \to F$ a morphism of Lie algebroids, and $\nabla$ a flat $F$-connection on $E$. Then we call the quadruple $(E, K, F, \nabla)$ a \textit{crossed module of Lie algebroids} if the following holds:
			\begin{align*}
				\nabla_X\mleft( \mleft[ \mu, \nu \mright]_E \mright)
				&=
				\mleft[ \nabla_X \mu, \nu \mright]_E
					+ \mleft[ \mu, \nabla_X\nu \mright]_E~,
				\\
				K\mleft( \nabla_X \mu \mright)
				&=
				\mleft[ X, K(\mu) \mright]_F~,
				\\
				\nabla_{K(\mu)} \nu 
				&=
				\mleft[ \mu, \nu \mright]_E
			\end{align*}
			for all $X \in \Gamma(F)$ and $\mu, \nu \in \Gamma(E)$.}
		\end{definitions}
		
		Observe, that $\basicc$ is an $E$-connection if $E$ is a BLA, and will thus now have no Leibniz rule. By definition we immediately get (keeping the same notation):
		
		\begin{corollaries}{Crossed modules as sandglass}{CrossedIsSand}
		{\color{Cred} \textit{To say that the quadruple $(E, K, F, \nabla)$ is a crossed module is equivalent to say that $\nabla$ is a flat Cartan $K$-connection whose basic connection $\basicc$ (both) is zero. In that case, choosing \color{Cred} $\zeta \equiv 0$ as primitive, there is a canonical sandglass sequence coupling $E$ and $F$ following $K$ for which the strict LAB $E_H$ of $E$ comes with the field of Lie brackets $H = - \mleft[ \cdot, \cdot\mright]_E$.}}
		\end{corollaries}
		
		\begin{proof}
		This is just a consequence of the definitions of $\basicc$ and $\basiccu$; then we choose a primitive $\zeta \equiv 0$ giving rise to a flat strict $K$-adjustment, inducing the canonical sandglass sequence provided in \ref{thm:CanonEx}. Henceforth, $H = t_{\basicc} = - \mleft[ \cdot, \cdot\mright]_E$.
		\end{proof}
		
		From the point of view provided in this paper, observe that we do not require flatness of $\nabla$, and we can allow $E$ to be a Lie algebroid itself. However, the mentioned references assume those stronger conditions in order to be able to define certain obstructions. As we pointed out, $\rmd^{\nabla^\zeta}\zeta$ obstructs sandglass sequences; we may explore in the next paper of this series whether or not this 3-form allows to generalise obstructions associated to crossed modules.
		
		Now, the fact that $H = - \mleft[ \cdot, \cdot\mright]_E$ and $\rho_E \equiv 0$ imply that $-K$ \textit{is} a morphism of Lie algebroids $E_H \to F$, while in general it is \textit{not}. Not only that: $\psi$, the natural projection to $F$ in the sandglass sequence \eqref{eq:TotalDiag}, is now also a morphism of Lie algebroids simply due to $\basicc = 0$ and $\zeta \equiv 0$ so that the whole sandglass sequence becomes a commuting diagram of Lie algebroid morphisms whose diagonals are short exact sequences. This begs the question: Does the sandglass in this case describe an equivalence? Indeed, it describes an equivalence of crossed modules, also called a \textit{butterfly} as introduced in \cite{Noohi:0910.1818, Aldrovandi:0808.3627}; see also \cite{Sheng:1003.1348}. Let us introduce this notion, but allowing $F$ to be a Lie algebroid as usual.
		
		\begin{definitions}{Butterflies}{Butterfliegen}
		\tcbsubtitle[before skip=\baselineskip]{\cite[the reference ``only'' considers Lie algebras]{Noohi:0910.1818}}
		\color{Cred}
		Let $(E, K, F, \nabla)$ and $(E', K', F', \nabla')$ be crossed modules of Lie algebroids over the same smooth manifold $M$. We speak of a \textit{butterfly (between $(E, K, F, \nabla)$ and $(E', K', F', \nabla')$)} if there is a Lie algebroid $A \to M$ and a commuting diagram of Lie algebroid morphisms
		\begin{equation*}
			\begin{tikzcd}[column sep=small]
			E \arrow[dd, "K", swap] \arrow[rd, "\hat\iota"] &  & E' \arrow[ld, "\iota", swap, hook] \arrow[dd, "K'"]
			\\
			& A \arrow[rd, "\psi", swap] \arrow[dl, "\scD", two heads] & 
			\\
			F & & F'
			\end{tikzcd}
		\end{equation*}
		such that
		\begin{align*}
		\mleft[ \xi, \hat\iota(\mu) \mright]_A
		&=
		\hat \iota \mleft( \nabla_{\scD(\xi)} \mu \mright)~,
		\\
		\mleft[ \xi, \iota(\nu) \mright]_A
		&=
		\iota \mleft( \nabla'_{\psi(\xi)} \nu \mright)
		\end{align*}
		for all $\mu \in \Gamma(E)$, $\nu \in \Gamma(E')$, and $\xi \in \Gamma(A)$.
		\end{definitions}
		
		Observe that the definition of butterfly only requires one diagonal to be short exact. Let us now prove that a sandglass is a special kind of butterfly in this case, making use of $\nabla^\zeta = \nabla$.
		
		\begin{lemmata}{Sandglasses as butterflies up to homotopy}{SandglassButterfly}
		\color{Cred}Let $(E, K, F, \nabla)$ be a crossed module of Lie algebroids. Then $(E_H, -K, F, \nabla)$ is also a crossed module of Lie algebroids where $E_H$ is the strict BLA of $E$, and the sandglass sequence of $(E, K, F, \nabla)$ as in \ref{cor:CrossedIsSand} is a butterfly between $(E, K, F, \nabla)$ and $(E_H, -K, F, \nabla)$ for which both diagonals are short exact sequences of Lie algebroids.
		\end{lemmata}
		
		\begin{proof}
		Using \ref{cor:CrossedIsSand}, we see that $-K\colon E_H \to F$ is a morphism of Lie algebroids because of $H = - \mleft[ \cdot, \cdot \mright]_E$ and $\rho_E \equiv 0$, and the canonical projection $\psi \colon A\coloneqq F\bowtie_{(\nabla, \zeta)} E \to F$ becomes a morphism of Lie algebroids, too, due to $\basicc = 0$ and $\zeta \equiv 0$. Observe that $(E_H, -K, F, \nabla)$ is obviously a crossed module of Lie algebroids by replacing $\mleft[ \cdot, \cdot \mright]_E$ with $H$ and $K$ with $-K$ in the equations of \ref{def:CrossedAlgoids}, making use of the fact that $(E, K, F, \nabla)$ is a crossed module. We have now a commuting diagram, by adjusting the sandglass sequence \eqref{eq:TotalDiag},
		\begin{equation*}
			\begin{tikzcd}[column sep=small]
			E \arrow[dd, "K", swap] \arrow[rd, "\hat\iota", hook] &  & E_H \arrow[ld, "\iota", swap, hook] \arrow[dd, "-K"]
			\\
			& A \arrow[rd, "\psi", swap, two heads] \arrow[dl, "\scD", two heads] & 
			\\
			F & & F
			\end{tikzcd}
		\end{equation*}
		where the arrows are defined as in the discussion around sandglass sequences and \ref{thm:CanonEx}. Making use of the definition of $A$, we can calculate that
		\begin{align*}
		\mleft[ (X, \mu), \hat\iota(\nu) \mright]_A
		&=
		\mleft[ (X, \mu), (0,\nu) \mright]_A
		=
		\mleft( 0 , \mleft[ \mu, \nu \mright]_E + \nabla_X\nu \mright)
		=
		\hat\iota\mleft( \nabla_{K(\mu)}\nu + \nabla_X\nu \mright)
		=
		\hat\iota\mleft( \nabla_{\scD(X, \mu)}\nu \mright)~,
		\\
		\mleft[ (X, \mu), \iota(\nu) \mright]_A
		&=
		\mleft[ (X, \mu), (-K(\nu),\nu) \mright]_A
		=
		\iota\mleft( \mleft[ \mu, \nu \mright]_E + \nabla_X\nu + \nabla_{K(\nu)} \mu \mright)
		=
		\iota\mleft( \nabla_X\nu \mright)
		=
		\iota\mleft( \nabla_{\psi(X, \mu)}\nu \mright)
		\end{align*}
		for all $\nu \in \Gamma(E)=\Gamma(E_H)$ and $(X, \mu) \in \Gamma(A)$, where we repeatedly used \ref{cor:CrossedIsSand}. This finishes the proof.
		\end{proof}
		
		We see, that in the context of crossed modules, sandglass sequences are equivalences of the crossed modules $(E, K, F, \nabla)$ and $(E_H, -K, F, \nabla)$. Thus, one may want to call $K$-twisted action Lie algebroids as in \ref{thm:ActionAlgebroidStructure} as \textit{twisted crossed modules of Lie algebroids}, or possibly also as \textit{crossed modules of Lie algebroids up to homotopy} as suggested by Iakovos Androulidakis in private communcation; and sandglass sequences as \textit{butterflies up to homotopy}.
		Indeed, let us conclude this paper by generalising the obstruction of a Lie group structure on a manifold $M$ (which is a special consequence of \textbf{Flat algebroid i)} of the introduction): The curvature of adjustments obstructs the structure of crossed modules. For this recall the definition of crossed modules of Lie groupoids: 
		
		\begin{definitions}{Crossed modules of Lie groupoids}{CrossedLieGroupoids}
		\tcbsubtitle[before skip=\baselineskip]{\cite[Definition 1.2, but we dropped the condition of $\Psi(\scG)$ being a closed submanifold because we will not discuss smoothness of maps defined on a quotient here]{androulidakis2005crossed}}
		\color{Cred}Let $\scG \to M$ be a Lie group bundle over a smooth manifold $M$, $\Omega$ a Lie groupoid over $M$, $\Psi\colon \scG \to \Omega$ a Lie groupoid morphism, and $\Phi\colon \Omega \to {{\rm Aut}}(\scG)$ an $\Omega$-representation on $\scG$, where ${{\rm Aut}}(\scG)$ is the Lie groupoid of automorphisms of $\scG$ covering diffeomorphisms of $M$. Then we call the quadruple $(\scG, \Psi, \Omega, \Phi)$ a \textit{crossed module of Lie groupoids} if the following holds:
			\begin{align*}
			\Psi(\Phi_h (g)) &= h ~ \Psi(g) ~ h^{-1}~,\\
			\Phi_{\Psi(g)}(g') &= gg'g^{-1} 
			\end{align*}
		for all $g, g' \in \scG$ over the same base point and $h \in \Omega$ whose source aligns with the base point of $g$. We will also say that this is a \textit{crossed module between $\scG$ and $\Omega$}.
		\end{definitions}
		
		This integrates \ref{def:CrossedAlgoids} and allows to integrate \ref{cor:CrossedIsSand}.
		
		\begin{lemmata}{Obstruction of crossed modules of Lie groupoids}{ObstrLieGroupoidCrossed}
		\color{Cred}Let $M$ be a smooth manifold, $\scG \to M$ a Lie group bundle whose fibres are connected and simply connected, and let $\Omega$ be a Lie groupoid over $M$ with connected and simply connected source\footnote{\color{Cred}We define the Lie algebroid of a Lie groupoid via the vertical structure of the source arrow.} fibres. Then there is a crossed module between $\scG$ and $\Omega$ if and only if there is a \color{Cred}Lie algebroid morphism $K\colon \frg \to F$ and a flat Cartan $K$-connection $\nabla$ with vanishing basic connection, where $\frg$ is the LAB of $\scG$ and $F$ the Lie algebroid of $\Omega$.
		\end{lemmata}
		
		\begin{proof}
		The proof will just consist of standard arguments of integrability in the $\Leftarrow$-direction. For this we will use Lie's Second Theorem, which also holds on the level of Lie groupoids; see for example \cite{Mackenzie:9712012, Moerdijk:0006042}. First of all, by \ref{cor:CrossedIsSand} we know that $(\frg, K, F, \nabla)$ is a crossed module of Lie algebroids. Now we integrate in the following way:
		\begin{itemize}
			\item View $\nabla$ as a morphism of Lie algebroids $F \to \scD_{{\rm Der}}(\frg)$ with values in the Lie algebroid of Lie bracket derivations of $\frg$. Then there is a unique morphism of Lie groupoids $\Phi \colon \Omega \to {\rm Aut}(\scG)$ integrating $\nabla$.
			\item Similarly, there is a unique Lie groupoid morphism $\Psi$ integrating $K$.
			\item Making use of the uniqueness and integrating the last two equations in \ref{def:CrossedAlgoids} in an ordinary way now shows that $(\scG, \Psi, \Omega, \Phi)$ is a crossed module of Lie groupoids; for the unfamiliar reader: The needed adjoint formulas for the Lie bracket essentially also hold for Lie groupoids but requires some exaplanation via frameworks like bisections, see for example \cite[\S 3.7]{0521499283}.
		\end{itemize}
		The other direction just follows by differentiation and by \ref{cor:CrossedIsSand}.
		\end{proof}
		
		Given a smooth, compact, connected and simply connected manifold $M$, we have now a cooking recipe for a structure as crossed module between $M$ as a Lie group and another connected and simply connected Lie group $H$:
		\begin{enumerate}
			\item First, $M$ has a Lie group structure if and only if it admits a flat Cartan connection (here: the ``classical'' notion).
			\item Now proceed as in \ref{lem:ObstrLieGroupoidCrossed}, setting $\scG = G$ and $\Omega = H$.
		\end{enumerate}
		
		However, this requires quite a lot of data. In the context of this paper the following theorem's setup may be more natural to study, directly generalising the classical statement when $M$ admits a Lie group structure; however, we will see that we have to ``twist'' the classical approach to avoid discussing parallel frames of Lie algebroid connections. See also \ref{rem:ParallelFrames} after the proof.
		
		\begin{theorems}{Obstruction of certain crossed modules of Lie groups}{ObstrCrossedModules}
		\color{Cred}Let $M$ be a smooth, compact, connected and simply connected manifold, $F\to M$ a Lie algebroid, and $K\colon \rmT M \to F$ a morphism of Lie algebroids. 
		
		If $\rmT M$ admits a flat Cartan $K$-connection $\nabla$, then $M$ is diffeomorphic to a Lie group $G$ and there is a crossed module of connected and simply connected Lie groups $(G, \Psi, H, \Phi)$ as in \ref{def:CrossedLieGroupoids} for which $G$ integrates a Lie algebra induced by a parallel frame of $\basicc$ on $\rmT M$, $H$ integrates a Lie algebra induced by a parallel frame of $\basicc$ on $F$, $K$ and $\nabla$ restrict to these Lie algebras, and $\Psi$ and $\Phi$ are integrals of $K$ and $\nabla$, respectively.

		Furthermore, if we additionally have $\nabla_\xi = 0$ for all $\xi\in \Gamma(F)$ with $\rho_F(\xi)=0$, then $H$ sits in a short exact sequence of Lie groups
		\begin{equation*}
			\begin{tikzcd}
				I \arrow[r, hook]& H \arrow[r, two heads] & G ~ ,
			\end{tikzcd}
		\end{equation*}
		where $I$ is the connected and simply connected Lie group integrating the structural isotropy Lie algebra\footnote{\color{Cred}As defined by the kernel of $\rho_F$.} of $F$; moreover, $\Psi$ splits this sequence inducing an isomorphism of Lie groups $H \cong G \ltimes I$ (semidirect product).
		\end{theorems}
		
		\begin{remarks}{Transitivity}{TransOfF}
		\color{Cred}That $K$ is a morphism implies that $K$ is a right-inverse of $\rho_F$, thus, $K$ is injective and $F$ transitive such that the structural isotropy Lie algebra is defined.
		\end{remarks}
		
		
		\begin{proof}[Proof of \ref{thm:ObstrCrossedModules}]

		\textbf{Step I: The Lie group structure $G$ on $M$.} Given a flat Cartan $K$-connection $\nabla$ we choose the trivial primitive $\zeta \equiv 0$ for it such that $(\nabla, \zeta \equiv 0)$ becomes a strict $K$-adjustment, and by \ref{prop:RevisitObstrActAlgoid} $t_{\basicc}$ with $\nabla t_{\basicc} = 0$ gives rise to a BLA structure on $\rmT M$ which is in fact an LAB (\ref{lem:BLALAB}). By the fact that $t_{\basicc} = - t_{\nabla_K}$ there is also an LAB structure associated to $t_{\nabla_K}$.
		
		Due to the vanishing basic curvature $\basiccu = 0$ we know that both $\basicc$ are flat, in particular we have a flat $\rmT M$-connection $\basicc$ on $\rmT M$, and due to the simply connectedness we have a global parallel frame $(e_i)_{i=1, \dots, \mathrm{dim}(M)}$ with $\basicc e_i = 0$, trivialising $\rmT M \cong M \times \frg$ as vector bundle for some vector space $\frg$ so that $(e_i)_i$ is the canonical frame of constant sections. Recall the first argument in \ref{rem:PossibleChoiceOfZetaKK}, $\basicc$ (on $\rmT M$) is in this case a flat Cartan connection on $\rmT M$, inducing the Lie group structure $G$. In more detail: By \ref{thm:ActionAlgebroidStructure}, $\rmT M$ is a $K$-twisted action Lie algebroid whose strict LAB $\rmT M_H$ is indeed w.r.t.\ $H = t_{\basicc}$, but observe that we can rewrite the $K$-twisted action Lie algebroid structure as in
		\begin{equation*}
		[\mu, \nu]_{\rmT M} 
		=
		t_{\basicc}(\mu, \nu) + \nabla_{K(\mu)} \nu - \nabla_{K(\nu)}\mu
		=
		t_{\nabla_K}(\mu, \nu) + \basicc_{\mu} \nu - \basicc_{\nu}\mu
		\end{equation*}
		for all $\mu ,\nu \in \frX(M)$, so that
		\begin{equation}\label{eq:TMaLAB}
		[e_i, e_j]_{\rmT M} 
		=
		t_{\nabla_K}(e_i, e_j)
		=
		\sum_{i,j,k}C^k_{ij} e_k
		\end{equation}
		for subordinate LAB-structure functions $C^k_{ij} \in C^\infty(M)$. By \ref{prop:HConstantWRTthebasicConn} it also holds that $\basicc t_{\basicc} = - \basicc t_{\nabla_K}= 0$, in particular, since $M$ is connected, $C^k_{ij}$ are constant and $\frg$ is henceforth the structural Lie algebra of the LAB structure induced by $t_{\nabla_K}$ (because $(e_i)$ is not only a finite-dimensional subspace of $\frX(M)$ but also a Lie subalgebra such that its associated trivialisation is also a trivialisation of $\rmT M$ viewed as LAB $\rmT M_H$). With $G$ we now denote the unique connected and simply connected Lie group integrating $\frg$.
		
		Now, the isomorphism $\rmT M \cong M \times \frg$ followed by the canonical projection to $\frg$ defines a 1-form $\omega$ with values in $\frg$ which is Maurer-Cartan, that is,
		\begin{equation*}
		\rmd \omega 
			+ \frac{1}{2} \mleft[ \omega \stackrel{\wedge}{,} \omega \mright]_\frg
		\cong
		\rmd^{\basicc} \omega 
			+ \frac{1}{2} t_{\nabla_K}\mleft( \omega \stackrel{\wedge}{,} \omega \mright)
		=
		0~,
		\end{equation*}
		and this follows as in \cite{Blaom:0404313}, where one has to see $\frg$ as constant sections in $\rmT M$ on the right side of $\cong$. As in \cite{Blaom:0404313}, it follows now that $M$ is diffeomorphic to $G$ (for example by \cite[\S 3, Theorem 8.7]{SharpeLieGroupStructure}).
		
		\textbf{Step II.a: The other Lie group $H$.} Now, recall that there is also the basic connection $\basicc$ on $F$ which is also flat as a vector bundle connection. Thus, fix another global parallel frame $(f_\alpha)_{\alpha = 1, \dots, {{\rm rk}}(F)}$ with $\basicc f_\alpha = 0$ trivialising $F \cong M \times \frh$ for a vector space $\frh$. In fact, since $K$ is injective (\ref{rem:TransOfF}), we construct the frame $(f_\alpha)$ as in the following: Recall that $K \circ \basicc = \basicc \circ K$, in particular $(K(e_i))_i$ is a global subframe of $F$ parallel to $\basicc$. Thus, extend $(K(e_i))$ to a global parallel frame $(f_\alpha)$ of $F$, and observe parallelity of $f_\alpha$ and the definition of $\basicc$ implies
		\begin{equation*}
		K\mleft( \nabla_{f_\alpha} \mu \mright)
		=
		\mleft[ f_\alpha, K(\mu) \mright]_F
		\end{equation*}
		for all $\mu \in \frX(M)$, in particular we get
		\begin{align*}
			\basicc_\mu \mleft(\mleft[ f_\alpha, f_\beta \mright]_F\mright)
			&=
			\mleft[ K(\mu), \mleft[ f_\alpha, f_\beta \mright]_F \mright]_F
				+ K\mleft( \nabla_{\mleft[ f_\alpha, f_\beta \mright]_F} \mu \mright)
			\\
			&=
			\mleft[ K(\mu), \mleft[ f_\alpha, f_\beta \mright]_F \mright]_F
				+ K\mleft( \mleft( \nabla_{f_\alpha}\nabla_{f_\beta} - \nabla_{f_\beta} \nabla_{f_\alpha} \mright) \mu \mright)
			\\
			&=
			\mleft[ K(\mu), \mleft[ f_\alpha, f_\beta \mright]_F \mright]_F
				+ \mleft[ f_\alpha, \mleft[ f_\beta, K\mleft(  \mu \mright) \mright]_F\mright]_F
				- \mleft[ f_\beta, \mleft[ f_\alpha, K\mleft(  \mu \mright) \mright]_F\mright]_F
			\\
			&=
			0~,
		\end{align*}
		where we used the Jacobi identity of $\mleft[ \cdot, \cdot \mright]_F$ and the flatness of $\nabla$; that is, similar to above, the frame of constant sections is a finite-dimensional Lie subalgebra of $\Gamma(F)$ implying that $\frh$ is a Lie algebra associated to it. By construction, $K$ restricts to an injective Lie algebra morphism $\frg \to \frh$ (recall Equation \eqref{eq:TMaLAB}), and we define $H$ as the connected and simply connected Lie group integrating $\frh$.
		
		\textbf{Step II.b: Special structure on $H$ given the extra assumption on $\nabla$.} Due to the transitivity of $F$, the kernel of $\rho_F$ is an LAB which we denote by $L$ with structural Lie algebra $\fri$, and $K$ splits the sequence
		\begin{equation}\label{eq:LFTM}
			\begin{tikzcd}
				L \arrow[r, hook]& F \arrow[r,"\rho_F", two heads] & \rmT M ~ ,
			\end{tikzcd}
		\end{equation}
		in particular, $F \cong \rmT M \oplus L$ as anchored vector bundle. Assume for now that $\nabla_\xi \mu = 0$ for all $\xi\in \Gamma(L)$. In this case we will be able to prove that $\basicc$ not only preserves the $\rmT M$-factor (which is due to $K \circ \basicc = \basicc \circ K$) but also the factor in $L$; for this we will actually also make use of the fact that $\nabla$ has now a global parallel frame $(\eta_i)_{i=1, \dots, \mathrm{dim}(M)}$, too, which is possible by the extra assumption: Given a splitting $F \cong \rmT M \oplus L$, $\nabla$ is just an ordinary flat $\rmT M$-connection giving rise to $(\eta_i)$ because the splitting preserves the anchors and $[\mu, \nu]_F \cong \mleft([\mu, \nu]_{\rmT M},0\mright)$ for all $\mu,\nu \in \frX(M)$ (recall Equation \eqref{eq:MackAlgebroid}); in other words, the parallel frame is inherited by the flatness of $\nabla_K$. This implies,
		\begin{equation*}
		\rho_F\mleft( \basicc_{\eta_i} \xi\mright)
		=
		\rho_F\mleft( \mleft[ K(\eta_i), \xi \mright]_F + K\mleft(\nabla_{\xi} \eta_i\mright) \mright)
		=
		0
		\end{equation*}
		for all $\xi \in \Gamma(L)$, therefore $\basicc$ restricts to $L$ making use of the fact that $\mu \mapsto \basicc_{\mu}$ is tensorial and $(\eta_i)$ is a global frame. In this case we can construct $(f_\alpha)$ by extending $(K(e_i))$ by adding a parallel frame of $L$. As above, $(f_\alpha)$ describes a finite-dimensional Lie algebra $\frh$, but since $L$ as the kernel of $\rho_F$ is an ideal in $F$ the part of the parallel frame in $L$ describes an ideal of $\frh$ which is isomorphic to $\fri$. We can extend Equation \eqref{eq:LFTM} to sections and then restrict to constant sections which is possible due to $\rho_F \circ K = \mathrm{id}_{\rmT M}$, i.e.\ $\rho_F$ now sends constant sections to constant sections by this definition of $(f_\alpha)$, resulting into a short exact sequence of Lie algebras,
		\begin{equation}\label{eq:ShortExactLiealgebras}
			\begin{tikzcd}
				\fri \arrow[r, hook, "\iota"]& \frh \arrow[r,"\rho_F", two heads] & \frg ~ .
			\end{tikzcd}
		\end{equation}
		Again, by construction $K$ restricts to a morphism of Lie algebras $\frg \to \frh$, but now also splits this sequence. Denote by $I$ and $H$ the unique connected and simply connected Lie group integrating $\fri$ and $\frh$, respectively. It is now a well-known fact that this split short exact sequence of Lie algebras integrates to the proposed short exact sequence of the corresponding connected and simply connected Lie groups, and integration of $K$ gives a splitting $\Psi$, casting $H$ into the form of a semidirect product $G \ltimes I$. 
		%
		
		\textbf{Step III: The crossed module.} Now back to the general situation, \textbf{dropping} the assumption that $\nabla_\xi \mu = 0$ for all $\xi\in \Gamma(L)$. First of all, observe that w.r.t.\ the constant frame $(e_i)_i$ (that is, on $\frg$) $\basicc e_i = 0$ and Equation \eqref{eq:TMaLAB} imply
		\begin{equation*}
		\nabla_{K(\mu)} \nu 
		= \mleft[ \mu, \nu \mright]_{\rmT M}
		= \mleft[ \mu, \nu \mright]_{\frg}
		\end{equation*}
		for all constant $\mu, \nu \in \frX(M)$. In the same manner, w.r.t.\ constant sections, that is, w.r.t.\ $\frg$ and $\frh$ we have due to $\basicc f_\alpha = 0$,
		\begin{equation*}
		K\mleft( \nabla_X \mu \mright)
		=
		\mleft[ X, K(\mu) \mright]_F
		=
		\mleft[ X, K(\mu) \mright]_\frh
		\end{equation*}
		for all constant $\mu \in \frX(M)$ and constant $X \in \Gamma(F)$, because $K$ maps constant sections to constant sections. Last, again w.r.t.\ constant sections, $\basiccu = 0$ restricts to $\nabla \mleft[ \cdot, \cdot \mright]_{\rmT M}=0$, and
		\begin{align*}
		\basicc_\mu \nabla_X \nu
		&=
		\mleft[ \mu, \nabla_X \nu \mright]_{\rmT M}
			+ \nabla_{K\mleft( \nabla_X \nu \mright)} \mu
		\\
		&=
		\mleft[ \mu, \nabla_X \nu \mright]_{\rmT M}
			+ \nabla_{\mleft[ X, K(\nu) \mright]_F} \mu
		\\
		&=
		\mleft[ \mu, \nabla_X \nu \mright]_{\rmT M}
			+ \nabla_X \nabla_{K(\nu)} \mu
			- \nabla_{K(\nu)} \nabla_X \mu
		\\
		&=
		\mleft[ \mu, \nabla_X \nu \mright]_{\rmT M}
			+ \nabla_X \mleft( \mleft[ \nu, \mu \mright]_{\rmT M} \mright)
			- \mleft[ \nu, \nabla_X \mu \mright]_{\rmT M}
		\\
		&=
		0~,
		\end{align*}
		where we used the constancy of the sections w.r.t.\ $\basicc$, flatness of $\nabla$, and $\nabla \mleft[ \cdot, \cdot \mright]_{\rmT M} = 0$. That is, $\nabla$ maps constant sections to constant sections, and thus restricts to an $\frh$-connection on $\frg$ with $\nabla \mleft[ \cdot, \cdot \mright]_{\frg} = 0$. By using \ref{lem:ObstrLieGroupoidCrossed} the proof is finished because we already argued that $K$ restricts to a Lie algebra morphism $\frg \to \frh$ so that we have everything to conclude the proof.
		\end{proof}
		
		\begin{remarks}{What parallel frame to take?}{ParallelFrames}
		\color{Cred}In contrast to ``classical'' approaches we took a parallel frame of $\basicc$ (on $\rmT M$) instead of $\nabla$, which allowed us to avoid discussing parallel frames for algebroid connections; this is not in contradiction to ``classical'' approaches, recall the first argument in \ref{rem:PossibleChoiceOfZetaKK}: For the obstruction of a Lie group structure on $M$ one needs a flat Cartan connection $\nabla$, in particular $F = \rmT M$ and $K = \mathrm{id}_{\rmT M}$. This recovers a symmetry: $\mleft(\basicc\mright)^{{\rm bas}} = \nabla$, furthermore, $\basicc$ is flat for Cartan connections, and since one assumes that $\nabla$ is flat it is straightforward to check that in this case $\basicc$ itself is a flat Cartan connection (\cite{Blaom:0404313}). Since one has $\nabla t_{\basicc} = \basicc t_{\basicc} = 0$, the Lie algebra inherited from parallel frames of $\nabla$ and $\basicc$ align. That is, in the classical theorem it does actually not matter on which connection one focuses, as long as both are flat. Alternatively, one may be able to argue with $\nabla_K$ as strict ${\rm id}_{\rmT M}$-adjustment.
		\end{remarks}
		
		This concludes this paper. One can probably drop compactness by a suitable definition of completeness of the involved connections, and dropping simply connectedness likely induces a crossed module structure on the universal cover of $M$ so that one could view the fundamental groupoid of $M$ as crossed module over a groupoid integrating $F$. But this will not be investigated here.
		
		}
		
    \section*{Acknowledgments}
    Special thanks to Camille Laurent-Gengoux for a lot of guidance and advice, also to Ralf Meyer, Iakovos Androulidakis, and to Ioan M{\u a}rcu{\c t} for fruitful discussions and guidance, and to Mehran Jalali Farahani, Hungrok Kim and Christian Sämann for their tremendous help regarding curved theories. I also want to thank Siye Wu, Mark John David Hamilton, Alessandra Frabetti and Chenchang Zhu for their great help and support. Also big thanks to my wife Rachelle, my mother, father, Dennis, Gregor, Marco, Nico, Jakob, Kathi, Konstantin, Lukas, Locki, Gareth, Philipp, Ramona, Annerose, Michael, Maxim, and Anna for all their love and support. Last but not least: A big thanks to my referee who provided me with a lot of guidance. I appreciate you all!
    
    \section*{Data and License Management}
    
    No additional research data beyond the data presented and cited in this work are needed to validate the research findings in this work. For the purposes of open access, the author has applied the \href{https://creativecommons.org/licenses/by/4.0/}{Creative Commons Attribution 4.0 International (CC~BY~4.0)} license to any author-accepted manuscript version arising from this work.

    \bibliography{bigone}
    
    \bibliographystyle{latexeu2}
    
\end{document}